\documentclass[a4paper,11pt]{amsart}
\usepackage{lipsum}
\usepackage{listings,xcolor}
\usepackage[affil-it]{authblk}
\usepackage{graphicx}
\usepackage{listings}
\usepackage{amssymb,amsfonts}
\usepackage{amsmath}
\usepackage[all,arc]{xy}
\usepackage{enumerate}
\usepackage{mathrsfs}
\usepackage[latin1]{inputenc}
\usepackage{palatino}
\usepackage{color}
\usepackage{graphicx}
\usepackage{graphics}
\usepackage{tikz}
\usepackage{epsfig}
\usepackage{url}
\usepackage{tikz}
\usepackage{mathtools}
\usepackage{epsfig}
\usepackage{listings,xcolor}
\usepackage{listings}
\usepackage{hyperref}
\usepackage[affil-it]{authblk}
\usetikzlibrary{graphs}
\usetikzlibrary{matrix}
\usetikzlibrary{arrows}
\usetikzlibrary{chains}

\theoremstyle{definition}

\newtheorem{thm}{Theorem}[section]

\theoremstyle{definition}
\newtheorem{defn}[thm]{Definition}

\theoremstyle{remark}

\theoremstyle{remark}

\numberwithin{equation}{section}
\makeatletter
\g@addto@macro{\endabstract}{\@setabstract}
\newcommand{\authorfootnotes}{\renewcommand\thefootnote{\@fnsymbol\c@footnote}}%
\makeatother

\definecolor{checkcolor}{rgb}{0.95, 0.95, 0.95}
\newsavebox{\definitionbox}
		{\end{minipage}\end{lrbox}%
	\begin{center}{\colorbox{checkcolor}{\usebox{\definitionbox}}}%
	\end{center}}
\begin{document}
	\begin{center}
		\LARGE 
		Weak Faddeev-Takhtajan-Volkov algebras; Lattice $W_n$ algebras; \par \bigskip
		
		\normalsize
		\authorfootnotes
		Farrokh Razavinia\footnote{Farrokh Razavinia}\textsuperscript{1} 
		
		\textsuperscript{1}Department of discrete mathematics, \\ Moscow Institute of Physics and Technology (MIPT) \par
	\textsuperscript{1} f.razavinia@phystech.edu\\
		\today
		
	\footnote{	The author would like to thank
			Professor Yaroslav Pugai for his helpful discussion during the
			preparation for this paper.}\textsuperscript{1}
		\footnote{ subjclass[2010]: Primary 16D10, 17B37, 81R50; Secondary  20G42.}
		\footnote{keywords: Lattice W algebras, quantum groups, Feigin's homomorphisms, Mathematica.}
	\end{center}
	
	\begin{abstract}
		In this paper, we will start by looking through our project's historical general view and then we will try to construct a new Poisson bracket on our simplest example $sl_2$ and then we will try to give a universal construction based on our universal variables and then will try to construct lattice $W_2$ algebras which will play a key role in our other constructions on lattice $W_3$ algebras and finally we will try to find the only nontrivial dependent generator of our lattice $W_4$ algebras and so on for lattice $W_n$ algebras.
		
		And at the end of this paper we will have appendix A, which will contain some parts of the Mathematica coding which we have used and have made for to find our algebra structures.
	\end{abstract}
	
	\section{Introduction}
	There is an old problem which has been considered and introduced by Boris Feigin in 1992. It has been born in its new formulation; on quantum Gelfand-Kirillov conjecture; in a public talk at RIMS in 1992 based on the nilpotent part of $U_q(g)$ i.e. $U_q(\mathfrak{n} )$ for $g$ a simple Lie algebra. \\
	Now, this problem is known as ``Feigin's Conjecture''.\\
	In the mentioned talk, Feigin proposed the existence of a certain family of homomorphisms on the quantized enveloping algebra $U_q(g)$  which will led us to a deffinition of lattice $W-$algebras.\\ 
	These ``homomorphisms'' has been turned to a very useful tool for to study the fraction field of quantized enveloping algebras. \cite{6}\\
	There have been many attempts to construct lattice $W$-algebras in Feigin's sence, which ensures the simplicity of the construction process of lattice $W$-algebra; for example the best known articles in the subject has been written by Kazuhiro Hikami and Rei Inoue who tried to obtain the algebra structure by using lax operators and generalized R matrices. \cite{7} \cite{8}\\
	Or Alexander Belov and Alexander Antonov and Karen Chaltikian, who first tried to follow Feigin's construction but finaly they also solved part of the conjecture by getting help of lax operators, and it made it very difficult to follow their publication.\cite{9} \cite{10} \\
	But here, in this paper, we will proceed and will introduce the simplest way of constructing such kind of algebras by just employing Feigin's homomorphisms and screening operators by defining a Poisson bracket on our variables just based on our Cartan matrix. \cite{1} \cite{2} \\
	We have to note that in \cite{2}, Yaroslav Pugai has constructed lattice $W_3$ algebras already, but here we will introduce its weaker version based on our newly defined Poisson bracket, constructed just based on the Cartan matrix $A_n$, which will make our job easier and more elegant.\\
	For to do this, let us set $C$ an arbitrary symmetrizable Cartan matrix of rank $r$ and let $n = n_{+}$ be the standard maximal nilpotent sub-algebra of the Kac-Moody algebra associated with $C$.\\
	So $n$ is generated by elements $E_1,\ldots, E_r$ which are satisfying in Serre relations, \cite{11} Where $r$ stands for $\text{rank}(C)$.\\
	In \cite{1}, we proved that screening operators $S_{X_{i}^{ji} } =  \sum\limits_{\substack{j \in \mathbb{Z} \\ \text{for} ~ i  ~ \text{fixed} }}^n X_{i}^{ji}$; for $X_{i}^{ji}$ generators of the $q-$commutative ring $\mathbb{C}_q[X_{i}^{ji}] := \frac{\mathbb{C}[X_{i}^{ji}]}{\left\langle X_{i}^{ji} X_{k}^{jk} - q^{<\alpha_i , \alpha_j > } X_{k}^{jk} X_{i}^{ji} \right\rangle }$ and for $<\alpha_i , \alpha_j > = a_{ij}$ the $ij$'s components of our Cartan matrix $C$; are satisfying in quantum Serre relations $\text{ad}_q (X_i)^{1 - a_{ij}} (X_j)$ for adjoint action $\text{ad}_q (X_i) (X_j) = X_i X_j - q^{a_{ij}} X_j X_i$ and $X_i \in (U_q)_{\alpha}$, $X_j \in (U_q)_{\beta} $ \cite{5}, for $(U_q)_{\alpha} = \{ u \in U_q(g) | q^{\mathfrak{h}} u q^{- \mathfrak{h}} = q^{\alpha (\mathfrak{h})} u ~~~~ \text{for all} ~ \mathfrak{h} \in \overset{ \vee }{P}  \} $ and $U_q(g) = \underset{\alpha ~ in ~ Q }{\oplus} (U_q)_{\alpha} $, for $Q = \underset{i \in I}{\oplus} \mathbb{Z}_{\alpha_i}$ the root lattice and for $\overset{ \vee }{P} $ a free Abelian group of rank $2 |I| - \text{rank} C$ with $\mathbb{Z}-$basis $\{ h_i | i \in I \} \cup \{ d_s | s = 1 , \cdots , |I| - \text{rank} C \}  $ and $\mathfrak{h} = {\rm I\!F} \otimes_{\mathbb{Z}} \overset{ \vee }{P}$ be the ${\rm I\!F}-$linear space spanned by $\overset{ \vee }{P}$. \cite{5} $\overset{ \vee }{P}$ will be called dual weight lattice and $\mathfrak{h}$ the Cartan subalgebra. And ${\rm I\!F}$ will stand for our ground field.\cite{5}\\
	Here for our Cartan matrix $C$, the quantum Serre relation will be \\
	\hspace*{3cm}$\text{ad}_q(X_i)^{1 - (-1)} (X_j) = \text{ad}_{q}^{2} (X_i) (X_j) $\\
	\hspace*{5.92cm} $= X_{i}^{2} X_j - [2]_q X_i X_j X_i + X_j X_{i}^{2}$\\ 
	\hspace*{5.92cm} $= X_{i}^{2} X_j - (q + q^{-1}) X_i X_j X_i + X_j X_{i}^{2}$\\
	Where $[2]_q$ stands for quantum number $[2]_q = \frac{q^2 - q^{-2}}{q - q^{-1}}$.\\
	And again as what we had in \cite{1}, we can define \\
	\hspace*{2.6cm} $U_q(n) : = \big \langle S_{X_{i}^{ji}} , S_{X_{k}^{jk}}  \mid (\text{ad}_q(S_{X_{i}^{ji}}))^2(S_{X_{k}^{jk}}) = 0 \big \rangle $,\\
	and for $\mathbb{C}_q[X]$ the quantum polynomial ring in one variable and twisted tensor product $\bar{\otimes}$, we can define \\
	\hspace*{1.6cm} $U_q(n) \bar{\otimes} \mathbb{C}_q[X_{l}^{jl}] := \big \langle  S_{X_{i}^{ji}} , S_{X_{k}^{jk}} , X_{l}^{jl}  \mid (\text{ad}_q(S_{X_{i}^{ji}}))^2(S_{X_{k}^{jk}}) = 0 $\\
	\hspace*{4.6cm}$ , S_{X_{i}^{ji}} X_{l}^{jl} = q^2 X_{l}^{jl} S_{X_{i}^{ji}} , S_{X_{k}^{jk}} X_{l}^{jl} = q^{-1} X_{l}^{jl} S_{X_{k}^{jk}}  \big \rangle $\\
	such that we have the following embeding
	$$ U_q(n) \hookrightarrow U_q(n) \bar{\otimes} \mathbb{C}_q[X_{l}^{jl}] \hookrightarrow U_q(n) \bar{\otimes} \mathbb{C}_q[X_{l}^{jl}] \bar{\otimes} \mathbb{C}_q[X_{m}^{jm}]  $$
	where $\mathbb{C}_q[X_{l}^{jl}] \bar{\otimes} \mathbb{C}_q[X_{m}^{jm}] = \mathbb{C} \big \langle X_{l}^{jl} , X_{m}^{jm} \mid X_{l}^{jl}  X_{m}^{jm} = q^{a_{lm}} X_{m}^{jm} X_{l}^{jl} $.\cite{1} \\
	Which will ensure the well-definedness of our definition of lattice $W-$algebras.
	\section{Weak Faddeev-Takhtajan-Volkov algebras}
	As it has been mentioned already in \cite{1}, the main tools which we will use, are difference equations, screening operators, Feigin's homomorphisms, adjoint actions, partial differential equations, and Cartan matrices.\\
	We know that from an abstract view $g = sl_{m+1}$ is an algebra related to the Cartan matrix $(a_{ij})_{i,j}$, for
	$a_{ij} =
	\begin{cases}
	2 & \text{if } i = j\\
	-1 & \text{if } |i - j| = 1\\
	0 & \text{if } |i-j|>1
	
	\end{cases}
	$ and so for $sl_2$ it will consist of just one row and one column, i.e. we have $A_1 = (2)$ and let us denote by $C\langle X\rangle$ the skew polynomial ring on generators $X=(X_i)_i$ labeled by $i \in \{ - \infty , \cdots -1, 0 , 1 , \cdots , +\infty \}$ and the defining $q-$commutation relations $X_i X_j = q^2 X_j X_i  ~~ \text{ for if} ~~ i \leq j$ with all having the same color.
	
	\begin{defn}
		Let's define our Poisson bracket as follows in the case of $sl_2$:\\
		\begin{equation}\label{Equ1}
			\begin{cases}
				\{X_i , X_j \} := 2 X_i X_j & \text{if } i < j\\
				\{ X_i , X_i \} := 0 
			\end{cases}
		\end{equation}
		
	\end{defn}

	The main problem is to find solutions of the system of difference equations from infinite number of non-commutative variables in quantum case and commutative variables in classical case. It is significant that commutation relations (\ref{Equ1}) depend just on the sign of the difference $(i-j)$ and is based on our Cartan matrix. We should try to find all solutions of the system:
	\begin{equation}\label{Equ2}
		\begin{cases}
			{\mathfrak{D}}_{x}^{(n)} \triangleleft \tau_1 = 0\\
			{H}_{x}^{(n)} \triangleleft   \tau_1 = 0
		\end{cases}
	\end{equation}
	
	Let us define our system of variables as follows 
	$$\vdots \hspace*{1cm} ~ \vdots ~ \hspace*{1cm} ~ \vdots ~ \hspace*{1cm} ~ \vdots ~ \hspace*{1cm} ~ \vdots \hspace*{1cm} ~ \vdots $$
	$$\cdots \hspace*{0.42cm} ~ X_{1}^{(11)} ~ \hspace*{0.42cm}~ X_{1}^{(21)} ~\hspace*{0.42cm} ~ X_{1}^{(31)} ~\hspace*{0.42cm} ~ X_{1}^{(41)} ~\hspace*{0.42cm} \cdots $$
	$$\cdots \hspace*{0.42cm} ~ X_{2}^{(12)} ~\hspace*{0.42cm} ~ X_{2}^{(22)} ~ \hspace*{0.42cm}~ X_{2}^{(32)} ~\hspace*{0.42cm} ~ X_{2}^{(42)} ~\hspace*{0.42cm} \cdots $$
	$$\cdots \hspace*{0.42cm} ~ X_{3}^{(13)} ~ \hspace*{0.42cm}~ X_{3}^{(23)} ~\hspace*{0.42cm} ~ X_{3}^{(33)} ~ \hspace*{0.42cm}~ X_{3}^{(43)} ~\hspace*{0.42cm} \cdots $$
	$$\cdots \hspace*{0.42cm} ~ X_{4}^{(14)} ~ \hspace*{0.42cm}~ X_{4}^{(24)} ~\hspace*{0.42cm} ~ X_{4}^{(34)} ~\hspace*{0.42cm} ~ X_{4}^{(44)} ~\hspace*{0.42cm} \cdots $$
	$$\vdots \hspace*{1cm} ~ \vdots ~ \hspace*{1cm} ~ \vdots ~ \hspace*{1cm} ~ \vdots ~ \hspace*{1cm} ~ \vdots \hspace*{1cm} ~ \vdots $$

	And let us equip this system of variables with lexicographic ordering, i.e. $j_{k_m} i < j_{k_n} i$ if $j_{k_m} < j_{k_n} $ and $ji_{k_m}<ji_{k_n}$ if $i_{k_m}< i_{k_n}$. And we need this kind of ordering because we have different kind of set of variables with a proper coloring such that each set has its own color different from its neighbors.\\
	We have $\tau_1 := \tau_1[\cdots , X_{1}^{(11)}, X_{1}^{(21)}, X_{1}^{(31)}, \cdots ,  X_{2}^{(12)}, X_{2}^{(22)}, X_{2}^{(32)} , \cdots ]$, a multi-variable function depend on $\{X_{i}^{(ji)}\}$'s for $i,j \in \{-\infty ,\cdots , 1, \cdots , n , \cdots , + \infty \}$ and $\mathfrak{D}_{x}^{(n)} $ comes from 
	\begin{equation}\label{Equ3}
		\{S_{X_{i}^{ji}} , \tau_1 \}_p = S_{X_{i}^{ji}} \tau_1 - p^{\text{deg}\tau_1 <\alpha_i, \alpha_j>} \tau_1 S_{X_{i}^{ji}}
	\end{equation}
	
	where $<\alpha_i, \alpha_j> = a_{ij}$ is related to our Cartan matrix and $S_{X_{i}^{ji}}$  is a screening operator on one of our variable sets, i.e. $S_{X_{i}^{ji}} = \sum\limits_{j \in \mathbb{Z}}^{} X_{i}^{ji}$. Then we will obtain the whole set of solutions by using the following shift operator:
	$$\tau_2 = \tau_1[X_{1}^{(11)} \rightarrow X_{1}^{(21)}, X_{1}^{(21)} \rightarrow X_{1}^{(31)},\cdots ],$$
	\begin{equation}\label{Equ4}
		\tau_3 = \tau_2[X_{1}^{(21)} \rightarrow X_{1}^{(31)}, X_{1}^{(31)} \rightarrow X_{1}^{(41)}, \cdots ] 
	\end{equation}
	$$\hspace{-4.6cm} \vdots$$

	\vspace*{0.3cm}
	\begin{defn}
		Let us define our lattice W-algebra based on its generators according to \cite{2} \cite{1}.\\
		Generators of lattice W-algebra associated with simple Lie algebra $g$ constitute of the functional basis of the space of invariants
		\begin{equation}\label{Equ5}
			\tau_i := \text{Inv}_{U_q(n_{+})}(\mathbb{C}_q[X_{i}^{ji} | i \in \mathbb{Z}])
		\end{equation} 
		with additional requirements
		\begin{equation}\label{Equ6}
			H_{X_{i}^{ji}} (\tau_i) = 0 ~~~~~~\ ~~~~ \ \text{and} ~~~~~~\ ~~~~ \ D_{X_{i}^{ji}} (\tau_i) = 0
		\end{equation} 
		where $H_{X_{i}^{ji}} $ and $D_{X_{i}^{ji}}$ will be specified later.
	\end{defn}
	Equation (\ref{Equ4}) means that the generators have to satisfy in quantum Serre relations and the first equation in (\ref{Equ6}) means that they should have zero degree.\\
	Here in this paper we just will work on the case where $g = \text{sl}_n$ and we will use $\tau_{i}^{(n)}$ instead of $\tau_i$. Where $(n) $ stands for $n$ in $\text{sl}_n$.\\
	\subsection{Lattice $W_2$ algebra}
	Let us first consider the $sl_2$ case and to simplifying the notations, let us consider our set of variables as $X_i := X_{i}^{ji} $.\\
	And as it has shown in \cite{1}, it is enough just to work with $S_{X_{i}^{ji}} =: S_{X_i} = \sum\limits_{i=1}^{3} X_i $, because the other parts for $i >3$ and $i<1$ will tend to zero.\\
	By setting $q = e^{-\mathfrak{h}}$, for the Planck constant $\mathfrak{h}$, we will try to find generators of our lattice $W_2$-algebra, in the case of $\text{sl}_2$.\\
	\textbf{First step:}\\
	First let us try to find $D_{X}^{(2)}$. \\
	For to do this and for simplicity, we will set $\tau_1 := \tau_1[ \cdots , X_1, X_2, X_3, , \cdots]$. And as it has been defined already, we have\\
	\hspace*{4cm} $D_{X}^{(2)} := \{ S_{X_i} , \tau_1 \}$\\
	\hspace*{4.82cm} $ = \{ X_1 + X_2 + X_3 , \tau_1 \}$\\
	\hspace*{4.82cm} $ = \{ X_1 , \tau_1\} + \{ X_2 , \tau_1\} + \{ X_3 , \tau_1\} $\\
	\begin{equation}\label{Equ7}
		\hspace*{0.84cm} = (D_{X_1} + D_{X_2} + D_{X_3}) \tau_1
	\end{equation}
	Now for to understand what is (\ref{Equ7}), we note that $D_{X_i} = \{ X_i , \tau_1\}$ and also we note that our function $\tau_1[\cdots , X_1 , X_2 , X_3 , \cdots ]$ is a polynomial function consist of powers of $X_i$. What I mean is that, it is enough to find $D_{X_i}$ on just powers of $X_j$ for different values of $j \in \mathbb{Z}$.\\
	So \\
	\begin{equation}\label{Equ8}
		(\ref{Equ7}) = \sum\limits_{j}^{} (\{ X_1 , X_{j}^{n} \} + \{ X_2 , X_{j}^{n} \} + \{ X_3 , X_{j}^{n} \} )
	\end{equation}
	Where according to rules which has been pointed out in \cite{1}, we have \\
	\hspace*{4cm} $\{ X_1 , X_{j}^{n}\} = X_1 X_{j}^{n} - q^{2n} X_{j}^{n} X_1  $\\
	\hspace*{5.62cm} $ = \begin{cases} 0, & \mbox{if } j > 1 \\ (1 - q^{4n}) X_1 X_{j}^{n} , & \mbox{if } j < 1 \\ (1 - q^{2n}) X_1 X_{j}^{n} , & \mbox{if } j = 1 \end{cases} $\\
	Where by setting $q = e^{-\mathfrak{h}}$ and letting $\mathfrak{h} = 1$ at the end, we will have:\\
	\hspace*{1cm} First case: $j > 1$;\\
	\hspace*{4cm} $\{ X_1 , X_{j}^{n} \} = 0$;\\
	\hspace*{1cm} Second case: $j < 1$;\\
	\hspace*{4cm} $\{ X_1 , X_{j}^{n} \} = (1 - e^{-4n \mathfrak{h}}) X_1 , X_{j}^{n} $\\
	\hspace*{5.52cm} $ \sim (1 -(1-4n \mathfrak{h})) X_1 , X_{j}^{n} $\\
	\hspace*{5.52cm} $ = 4n \mathfrak{h} X_1 , X_{j}^{n} \sim 4n  X_1 , X_{j}^{n}$ \\
	\hspace*{5.52cm} $ = 4 X_1 X_j \frac{\partial X_{j}^{n}}{\partial X_j}$.\\
	\hspace*{1cm} Third case: $j = 1$;\\
	\hspace*{4cm} $\{ X_1 , X_{1}^{n} \} = (1 - q^{2n}) X_1 X_{1}^{n} $\\
	\hspace*{5.52cm} $ = (1 - e^{-2 n \mathfrak{h}}) X_1 X_{1}^{n} $\\
	\hspace*{5.52cm} $ \sim (1 - (1 -2 n \mathfrak{h})) X_1 X_{1}^{n} $\\
	\hspace*{5.52cm} $  = 2 n \mathfrak{h} X_1 X_{1}^{n} \sim 2n X_1 X_{1}^{n} $\\
	\hspace*{5.52cm} $ = 2 X_{1}^{2} \frac{\partial  X_{1}^{n} }{\partial X_1}$.\\
	And so we have\\
	\hspace*{4cm} $  (\ref{Equ8}) = \{ X_1 , X_{1}^{n} \} + \sum\limits_{j < 1}^{} \{ X_1 , X_{j}^{n} \} + \sum\limits_{j > 1}^{} \{ X_1 , X_{j}^{n} \}  $\\
	\hspace*{4.82cm} $  + \{ X_2 , X_{2}^{n} \} + \sum\limits_{j < 2}^{} \{ X_2 , X_{j}^{n} \} + \sum\limits_{j > 2}^{} \{ X_2 , X_{j}^{n} \}  $\\
	\hspace*{4.82cm} $  + \{ X_3 , X_{3}^{n} \} + \sum\limits_{j < 3}^{} \{ X_3 , X_{j}^{n} \} + \sum\limits_{j > 3}^{} \{ X_3 , X_{j}^{n} \}  $\\
	\hspace*{4.82cm} $  = 2 X_{1}^{2} \frac{\partial}{\partial X_1}  + 0 + 0 $\\
	\hspace*{4.82cm} $  + 2 X_{2}^{2} \frac{\partial}{\partial X_2}  + 4 X_2 X_1 \frac{\partial}{\partial X_1} + 0 $\\
	\hspace*{4.82cm} $  + 2 X_{3}^{2} \frac{\partial}{\partial X_3}  + 4 X_3 X_2 \frac{\partial}{\partial X_2} + 4 X_3 X_1 \frac{\partial}{\partial X_1}  $\\
	\hspace*{4.82cm} $  = 2X_1 (X_1 + 2 X_2 + 2 X_3) \frac{\partial}{\partial X_1} + 2 X_2 (X_2 + 2 X_3) \frac{\partial}{\partial X_2}  $\\
	\hspace*{4.82cm} $+2 X_{3}^{2} \frac{\partial}{\partial X_3} .$\\
	So we found $D_{X}^{(2)}$ which is as follows and we can omit 2, because finally we will make the action equal to zero and we can cancel 2  from both sides. So we have\\
	\begin{equation}\label{Equ9}
		D_{X}^{(2)} = X_1 (X_1 + 2 X_2 + 2 X_3) \frac{\partial}{\partial X_1} +  X_2 (X_2 + 2 X_3) \frac{\partial}{\partial X_2}  + X_{3}^{2} \frac{\partial}{\partial X_3}
	\end{equation}
	\textbf{Second step:}\\
	Now we will try to find $H_{X}^{(2)}$. \\
	For to find $H_{X}^{(2)}$, we note that it resembles the degree of our polynomial function. So if for example $H_{X}^{(2)}$ acts on $X_{1}^{n} X_{2}^{m} X_{3}^{l}$, then we should get $(n+m+l)$. \\
	So let us define: 
	\begin{equation}\label{Equ10}
		H_{X}^{(2)} := \sum\limits_{i}^{} X_i \frac{\partial}{\partial X_i}
	\end{equation}
	and then we have;\\
	\hspace*{2cm} $H_{X}^{(2)} (X_{1}^{n} X_{2}^{m} X_{3}^{l}) = (\sum_{i}^{} X_i \frac{\partial}{\partial_{X_i}}) (X_{1}^{n} X_{2}^{m} X_{3}^{l}) $\\
	\hspace*{4.68cm}$= \sum_{i}^{} X_i \frac{\partial X_{1}^{n} X_{2}^{m} X_{3}^{l} }{\partial_{X_i}} $ \\
	\hspace*{4.62cm} $= X_1  \frac{\partial X_{1}^{n} X_{2}^{m} X_{3}^{l} }{\partial_{X_1}} + X_2 \frac{\partial X_{1}^{n} X_{2}^{m} X_{3}^{l} }{\partial_{X_2}} + X_3 \frac{\partial X_{1}^{n} X_{2}^{m} X_{3}^{l} }{\partial_{X_3}} $\\
	\hspace*{4.66cm}$= n X_{1}^{n} X_{2}^{m} X_{3}^{l} + m X_{1}^{n} X_{2}^{m} X_{3}^{l}+ l X_{1}^{n} X_{2}^{m} X_{3}^{l} $\\
	\hspace*{4.69cm}$= (n+m+l) X_{1}^{n} X_{2}^{m} X_{3}^{l}$.\\
	Which gives us \\
	\hspace*{4cm} $H_{X}^{(2)} (X_{1}^{n} X_{2}^{m} X_{3}^{l}) = (n+m+l) X_{1}^{n} X_{2}^{m} X_{3}^{l}$\\
	and on the other side we have\\
	\hspace*{0.2cm} $ (n+m+l) X_{1}^{n} X_{2}^{m} X_{3}^{l} = n X_1 X_{1}^{n-1} X_{2}^{m} X_{3}^{l}  + m X_{1}^{n} X_2 X_{2}^{m-1} X_{3}^{l} + l X_{1}^{n} X_{2}^{m} x_3  X_{3}^{l-1}  $\\
	\hspace*{3.48cm} $= X_1 \frac{X_{2}^{m} X_{3}^{l} \partial X_{1}^{n} }{\partial_{X_1}} + X_2 \frac{X_{1}^{n} X_{3}^{l} \partial  X_{2}^{m} }{\partial_{X_2}} + X_3 \frac{X_{1}^{n} X_{2}^{m} \partial  X_{3}^{l} }{\partial_{X_3}}  $\\
	\hspace*{3.5cm}$=  X_1 \frac{ \partial  }{\partial_{X_1}} + X_2 \frac{\partial  }{\partial_{X_2}} + X_3 \frac{\partial  }{\partial_{X_3}}$\\
	Which gives us \\
	\hspace*{4cm} $(n+m+l) X_{1}^{n} X_{2}^{m} X_{3}^{l}  = \sum\limits_{i}^{} X_i \frac{\partial}{\partial X_i} $.\\
	And it shows that $(\ref{Equ10})$ is well defined.\\
	Now the only thing that remains is just to find the solutions of the following system of 2-linear homogeneous equations in one unknown $\tau_1$:\\
	\begin{equation}\label{Equ11}
		\begin{cases} (X_1 ( X_1 + 2 X_2 + 2 X_3) \frac{\partial }{\partial X_1} +\hspace{-0.1cm} X_2 ( X_2 + 2 X_3) \frac{\partial }{\partial X_2} +  X_{3}^{2} \frac{\partial }{\partial X_3} ) \tau_1[\cdots \hspace{-0.1cm}, X_1, X_2,  X_3, \\ \hspace{0.4cm} \cdots ] \hspace{-0.1cm}=0, & \mbox{ }  \\ ( X_1 \frac{ \partial  }{\partial_{X_1}} + X_2 \frac{\partial  }{\partial_{X_2}} + X_3 \frac{\partial  }{\partial_{X_3}}) \tau_1[\cdots , X_1, X_2,  X_3,  \cdots ] =0. & \mbox{ }  \end{cases}
	\end{equation}
	Now the goal is to find such $\tau_1[\cdots , X_1, X_2,  X_3,  \cdots ]$ which satisfies in our system of equations (\ref{Equ11}).\\
	The second equation ensures that the solution has degree 0 and also the partial differentials will give us a multi-variable function dependent on just $X_1, X_2,  X_3$.\\
	The system of PDEs (\ref{Equ11}) can be solved using the procedure described in Chapter V, Section IV of \cite{3}. \\
	And after doing some calculation in mathematica it become clear that the system (\ref{Equ11}) has only one functional dependent nontrivial solution:
	\begin{equation}\label{Equ12}
		\tau_{1}^{(2)}[X_1, X_2, X_3] = \frac{(X1 + X2 )(X2 + X3)}{X2 (X1 + X2 + X3)} = \frac{(\sum_{\substack{1 \leq i_1 \leq 2}} X_{i_1}^{(1)})(\sum_{\substack{1 \leq i_1 \leq 2}} X_{i_1 + 1}^{(1)})}{X_{2}^{(1)} (\sum_{\substack{1 \leq i_1 \leq 3}} X_{i_1}^{(1)} )}. 
	\end{equation}
	And again as before, $(2)$ goes back to 2 in $Sl_2$ and  $1$ is a default index which will be used later for to employ shifting operator.\\
	According to the number of variables, we will have two shifts and then everything will be in a loop.\\
	So here in $sl_2$ case we have three solutions for our system of linear equations $(\ref{Equ11})$ which  belong to the fraction ring of polynomial functions:
	\begin{equation}\label{Equ13}
		\begin{cases}
			\tau_{1}^{(2)}[X_1, X_2, X_3]= \frac{(\sum_{\substack{1 \leq i_1 \leq 2}} X_{i_1}^{(1)})(\sum_{\substack{1 \leq i_1 \leq 2}} X_{i_1 + 1}^{(1)})}{X_{2}^{(1)} (\sum_{\substack{1 \leq i_1 \leq 3}} X_{i_1}^{(1)} )} ; \\
			\tau_{2}^{(2)}[X_2, X_3, X_4]= \frac{(\sum_{\substack{2 \leq i_1 \leq 3}} X_{i_1}^{(1)})(\sum_{\substack{2 \leq i_1 \leq 3}} X_{i_1 + 1}^{(1)})}{X_{2}^{(1)} (\sum_{\substack{2 \leq i_1 \leq 4}} X_{i_1}^{(1)} )} ; \\
			\tau_{3}^{(2)}[X_3, X_4, X_5]= \frac{(\sum_{\substack{3 \leq i_1 \leq 4}} X_{i_1}^{(1)})(\sum_{\substack{3 \leq i_1 \leq 4}} X_{i_1 + 1}^{(1)})}{X_{2}^{(1)} (\sum_{\substack{3 \leq i_1 \leq 5}} X_{i_1}^{(1)} )} . 
		\end{cases}
	\end{equation}
	
	We go to define our non-commutative Poisson algebra according to definition of Poisson brackets given by Poisson himself \cite{4} with the difference that here we work on $q-$commutative ring $\frac{\mathbb{C}[X_{i}^{ji}]}{X_{i}^{ji} X_{k}^{jk} - q^{<\alpha_i , \alpha_k>} X_{k}^{jk} X_{i}^{ji}}$,  based on the generators which are the solutions of PDEs system  $(\ref{Equ2})$. \\
	For to do this we will use the following bracket:
	\begin{equation}\label{Equ14}
		F_{j}^{(n)} := \{ \tau_{i}^{(n)} , \tau_{j}^{(n)} \} = \sum\limits_{i}^{} \frac{\partial \tau_{i}^{(n)} }{\partial X_{i}}  \sum\limits_{j}^{} \frac{\partial \tau_{j}^{(n)}}{\partial X_{j}}  \{ X_{i} , X_{j}  \},
	\end{equation}
	where $\{ X_{i} , X_{j}  \}$ is our previously defined Poisson bracket on our set of variables.\\
	For instance in the case of $sl_2$  we have\\
	$$\{ \tau_{1}^{(2)} , \tau_{2}^{(2)} \} =\Big( \frac{\partial \tau_{1}^{(2)}}{\partial X_{1}}   \Big) \Big( \frac{\partial \tau_{2}^{(2)}}{\partial X_{2}} \{ X_1 , X_2 \}  + \frac{\partial \tau_{2}^{(2)}}{\partial X_{3}} \{ X_1 , X_3 \}  + \frac{\partial \tau_{2}^{(2)}}{\partial X_{2}} \{ X_1 , X_4 \}   \Big)  $$
	$$\hspace*{1.7cm}+  \Big(  \frac{\partial \tau_{1}^{(2)}}{\partial X_{2}}  \Big) \Big( \frac{\partial \tau_{2}^{(2)}}{\partial X_{2}} \{ X_2 , X_2 \}  + \frac{\partial \tau_{2}^{(2)}}{\partial X_{3}} \{ X_2 , X_3 \}  + \frac{\partial \tau_{2}^{(2)}}{\partial X_{2}} \{ X_2 , X_4 \}   \Big) $$
	$$\hspace*{1.7cm} + \Big(  \frac{\partial \tau_{1}^{(2)}}{\partial X_{3}}  \Big) \Big( \frac{\partial \tau_{2}^{(2)}}{\partial X_{2}} \{ X_3 , X_2 \}  + \frac{\partial \tau_{2}^{(2)}}{\partial X_{3}} \{ X_3 , X_3 \}  + \frac{\partial \tau_{2}^{(2)}}{\partial X_{2}} \{ X_3 , X_4 \}   \Big) $$
	$$\hspace*{1.7cm}  = \Big( \frac{\partial \tau_{1}^{(2)}}{\partial X_{1}}   \Big) \Big( \frac{\partial \tau_{2}^{(2)}}{\partial X_{2}} (2 X_1 X_2)  + \frac{\partial \tau_{2}^{(2)}}{\partial X_{3}} (2 X_1 X_3)  + \frac{\partial \tau_{2}^{(2)}}{\partial X_{2}} (2 X_1 X_4)   \Big) $$
	$$ \hspace*{0.6cm}+  \Big(  \frac{\partial \tau_{1}^{(2)}}{\partial X_{2}}  \Big) \Big( \frac{\partial \tau_{2}^{(2)}}{\partial X_{2}} (0)  + \frac{\partial \tau_{2}^{(2)}}{\partial X_{3}} (2X_2 X_3)  + \frac{\partial \tau_{2}^{(2)}}{\partial X_{2}} (2 X_2 X_4)   \Big) $$
	$$\hspace*{0.9cm} + \Big(  \frac{\partial \tau_{1}^{(2)}}{\partial X_{3}}  \Big) \Big( \frac{\partial \tau_{2}^{(2)}}{\partial X_{2}} (-2 X_3 X_2)  + \frac{\partial \tau_{2}^{(2)}}{\partial X_{3}} (0) + \frac{\partial \tau_{2}^{(2)}}{\partial X_{2}} (2 X_3 X_4)   \Big)$$
	$$\hspace*{-2.0cm} = 2 \frac{ X_1 X_{2}^{2} X_{3}^{2} X_4 (X_1 + X_2 + X_3 + X_4)}{(X_1 + X_2)^2 (X_2 + X_3)^3 (X_3 + X_4)^2}$$
	So we have
	\begin{equation}\label{Equ15}
		F_{2}^{(2)} = \{ \tau_{1}^{(2)} , \tau_{2}^{(2)} \} = \frac{2 X_1 X_{2}^{2} X_{3}^{2} X_4 (X_1 + X_2 + X_3 + X_4)}{(X_1 + X_2)^2 (X_2 + X_3)^3 (X_3 + X_4)^2}
	\end{equation}
	And it is enough to find our brackets just based on the first generator, because after that we are able to find other brackets based on the other generators, so for $\tau_{3}^{(2)}$ in an almost same process  we have:\\
	\hspace*{1.14cm}$F_{3}^{(2)} = \{ \tau_{1}^{(2)} , \tau_{3}^{(2)} \} $\\
	$$\hspace*{1.6cm}= \Big( \frac{\partial \tau_{1}^{(2)}}{\partial X_{1}}   \Big) \Big( \frac{\partial \tau_{3}^{(2)}}{\partial X_{3}} \{ X_1 , X_3 \}  + \frac{\partial \tau_{3}^{(2)}}{\partial X_{4}} \{ X_1 , X_4 \}  + \frac{\partial \tau_{3}^{(2)}}{\partial X_{5}} \{ X_1 , X_5 \}   \Big)  $$
	$$ \hspace*{1.6cm} + \Big( \frac{\partial \tau_{1}^{(2)}}{\partial X_{2}}   \Big) \Big( \frac{\partial \tau_{3}^{(2)}}{\partial X_{3}} \{ X_2 , X_3 \}  + \frac{\partial \tau_{3}^{(2)}}{\partial X_{4}} \{ X_2 , X_4 \}  + \frac{\partial \tau_{3}^{(2)}}{\partial X_{5}} \{ X_2 , X_5 \}   \Big)  $$
	$$ \hspace*{1.6cm} + \Big( \frac{\partial \tau_{1}^{(2)}}{\partial X_{3}}   \Big) \Big( \frac{\partial \tau_{3}^{(2)}}{\partial X_{3}} \{ X_3 , X_3 \}  + \frac{\partial \tau_{3}^{(2)}}{\partial X_{4}} \{ X_3 , X_4 \}  + \frac{\partial \tau_{3}^{(2)}}{\partial X_{5}} \{ X_3 , X_5 \}   \Big)  $$ 
	$$\hspace*{1.6cm}  = \Big( \frac{\partial \tau_{1}^{(2)}}{\partial X_{1}}   \Big) \Big( \frac{\partial \tau_{3}^{(2)}}{\partial X_{3}} (2 X_1 X_3)  + \frac{\partial \tau_{3}^{(2)}}{\partial X_{4}} (2 X_1 X_4)  + \frac{\partial \tau_{3}^{(2)}}{\partial X_{5}} (2 X_1 X_5)   \Big) $$
	$$ \hspace*{1.4cm} + \Big( \frac{\partial \tau_{1}^{(2)}}{\partial X_{2}}   \Big) \Big( \frac{\partial \tau_{3}^{(2)}}{\partial X_{3}} (2 X_2 X_3)  + \frac{\partial \tau_{3}^{(2)}}{\partial X_{4}} (2 X_2 X_4) + \frac{\partial \tau_{3}^{(2)}}{\partial X_{5}} (2 X_2 X_5)   \Big)  $$
	$$ \hspace*{0.5cm} + \Big( \frac{\partial \tau_{1}^{(2)}}{\partial X_{3}}   \Big) \Big( \frac{\partial \tau_{3}^{(2)}}{\partial X_{3}} (0)  + \frac{\partial \tau_{3}^{(2)}}{\partial X_{4}} (2 X_3 X_4)  + \frac{\partial \tau_{3}^{(2)}}{\partial X_{5}} (2 X_3 X_5)  \Big)  $$
	\begin{equation}\label{Equ16}
		\hspace*{-1.0cm} =  \frac{ -2 X_1 X_2 X_{3}^{2} X_4 X_5}{(X_1 + X_2) (X_2 + X_3)^2 (X_3 + X_4)^2 (X_4 + X_5)}.
	\end{equation}
	We have to note that we almost are done with our Poisson algebra in $\text{sl}_2$ case, but for the further plan i.e. to find our Volterra system, the differential-difference chain of non-linear equations
	\begin{equation}\label{Equ17}
		\begin{cases}
			H = \sum\limits_{i}^{} [\ln(\tau_i)]; \\
			\dot{\tau_j} = \{ \tau_j , H \} = \tau_j \times \sum\limits_{i}^{} \Gamma_i;
		\end{cases}
	\end{equation}
	where $\Gamma_i$ stands for $\frac{{\tau_1 , \tau_i}}{\tau_1 \tau_i}$ \cite{2}, we have to write down the brackets $\{ \tau_1 , \tau_i\} $ in terms of their decompositions to $\tau_j$'s for $1 \leq j \leq i$.\\
	So we need to write it as the decomposition of our generators and it will be done by using the Mathematica coding which we have produced in Appendix A. 
	
	And the result is as follows:
	\begin{equation}\label{Equ18}
		\begin{cases}
			F_{2}^{(2)} = \{ \tau_{1}^{(2)} , \tau_{2}^{(2)} \} = 2 (1 - \tau_{1}^{(2)} ) (1 - \tau_{2}^{(2)}) (-1 + \tau_{1}^{(2)}  +  
			\tau_{2}^{(2)}); \\
			F_{3}^{(2)} = \{ \tau_{1}^{(2)} , \tau_{3}^{(2)} \} = -2 (1 - \tau_{1}^{(2)} ) (1 - \tau_{2}^{(2)}) (1 - \tau_{3}^{(2)}); \\
			F_{i}^{(2)} = \{ \tau_{1}^{(2)} , \tau_{i}^{(2)} \} = 0 & \hspace*{-1cm} \text{for } |i - 1| \geq 3; 
		\end{cases}
	\end{equation}
	
	This result is weaker than the Faddeev-Takhtajan-Volkov algebra which has been mentioned in \cite{2} and if we continue this for $sl_3$, then we will have again a weaker version of what which has been mentioned in \cite{2}.
	\subsection{Lattice $W_3$ algebra}
	In this case we will use the following defined Poisson bracket based on Cartan matrix $A_2 =  \left[ \begin{array}{cc}
	$2$&$-1$ \\
	$-1$&$2$ \\
	\end{array} \right],$ but for to do this according to our previous ordering and list of variables, let us for simplicity set our variables as follows\\
	Set $X_{i}^{(1i)} := X_i$ and $X_{i}^{(2i)} := Y_i$.\\
	\begin{defn}
		Let's define our Poisson bracket as follows in the case of $sl_3$:
		\begin{equation}\label{Equ19}
			\begin{cases}
				\{X_i , X_j \} := 2 X_i X_j & \text{if } i < j;\\
				\{Y_i , Y_j \} := 2 Y_i Y_j & \text{if } i < j;\\
				\{ X_i , X_i \} := 0; \\
				\{ Y_i , Y_i \} := 0; \\
				\{ X_i , Y_j \} := X_i Y_j & \text{if } i > j;\\
				\{ X_i , Y_j \} := - X_i Y_j & \text{if } i \leq j;
			\end{cases}
		\end{equation}

	\end{defn}
	

	And instead of $(\ref{Equ1})$ we will have the following $q-$commutation relations 
	\begin{equation}\label{Equ20}
		\begin{cases}
			X_i X_j = q^2 X_j X_i  & \text{if } i \leq j;\\
			Y_i Y_j = q^2 Y_j Y_i  & \text{if } i \leq j;\\
			X_i Y_j = q^{-1} Y_j X_i & \text{if } i \leq j.
		\end{cases}
	\end{equation}

	And we will get the following equations in a same manner as in $\text{sl}_2$:\\
	\hspace*{1cm} First case: $i < j$;\\
	\hspace*{4cm} $ \{ X_i , Y_{j}^{n}\} = X_i Y_{j}^{n} - q^{-n} Y_{j}^{n} X_i $\\
	\hspace*{5.82cm} $=  X_i Y_{j}^{n} - q^0 X_i X_{j}^{n} $ \\
	\hspace*{5.82cm} $ = 0 $\\
	\hspace*{1cm} Second case: $i \geq j$;\\
	
	\hspace*{4cm} $ \{ X_i , Y_{j}^{n}\} = X_i Y_{j}^{n} - q^{-n} Y_{j}^{n} X_i $\\
	\hspace*{5.85cm} $=  (1 - q^{-2n} ) X_i Y_{j}^{n}  $\\
	\hspace*{5.85cm} $= (1 - e^{2n \mathfrak{h}} ) X_i Y_{j}^{n}  $ \\
	\hspace*{5.85cm} $ \sim (1 - (1 + 2n \mathfrak{h} )) X_i Y_{j}^{n}  $\\
	\hspace*{5.85cm} $= -2n \mathfrak{h} X_i Y_{j}^{n}  $ \\
	\hspace*{5.85cm} $\sim - 2n X_i Y_{j}^{n} $\\
	\hspace*{5.85cm} $ = -2 X_i Y_j \frac{\partial Y_{j}^{n}}{\partial Y_{j}}$ \\
	
	\begin{equation}\label{Equ21}
		\begin{cases}
			\{ X_i , X_{j}^{n}\} = 0 & \text{if } i \leq j;\\
			\{ X_i , X_{j}^{n}\} = 4 X_i X_j \frac{\partial X_{j}^{n}}{\partial X_{j}} & \text{if } i > j;\\
			\{ X_i , Y_{j}^{n}\} = 0 & \text{if } i < j;\\
			\{ X_i , Y_{j}^{n} \} = -2 X_i Y_j \frac{\partial Y_{j}^{n}}{\partial Y_{j}} & \text{if } i \geq j; \\
			\{ Y_j , X_{i}^{n} \} = -2 Y_j X_i \frac{\partial X_{i}^{n}}{\partial X_{i}} & \text{if } i \leq j;
		\end{cases}
	\end{equation}

	As in $sl_2$ case we will try to find $H_{X}^{(3)}$ as follows:\\
	$\{ X_{1}^{\alpha_1} X_{2}^{\alpha_2} X_{3}^{\alpha_3}  Y_{1}^{\beta_1} Y_{2}^{\beta_2} Y_{3}^{\beta_3}  , X_0 \} $\\
	$= X_{1}^{\alpha_1} X_{2}^{\alpha_2} X_{3}^{\alpha_3}  Y_{1}^{\beta_1} Y_{2}^{\beta_2} Y_{3}^{\beta_3} X_0 - X_0 X_{1}^{\alpha_1} X_{2}^{\alpha_2} X_{3}^{\alpha_3}  Y_{1}^{\beta_1} Y_{2}^{\beta_2} Y_{3}^{\beta_3}$\\
	$= (1 - q^{2 \alpha_1 + 2 \alpha_2 + 2 \alpha_3 - \beta_1 - \beta_2 - \beta_3  })  X_{1}^{\alpha_1} X_{2}^{\alpha_2} X_{3}^{\alpha_3}  Y_{1}^{\beta_1} Y_{2}^{\beta_2} Y_{3}^{\beta_3} X_0 $\\
	$\sim (1 - (1 - n \mathfrak{h}(2 \alpha_1 + 2 \alpha_2 + 2 \alpha_3 - \beta_1 - \beta_2 - \beta_3)) ) X_{1}^{\alpha_1} X_{2}^{\alpha_2} X_{3}^{\alpha_3}  Y_{1}^{\beta_1} Y_{2}^{\beta_2} Y_{3}^{\beta_3} X_0 $\\
	$= (2 \alpha_1 + 2 \alpha_2 + 2 \alpha_3 - \beta_1 - \beta_2 - \beta_3) n \mathfrak{h} X_{1}^{\alpha_1} X_{2}^{\alpha_2} X_{3}^{\alpha_3}  Y_{1}^{\beta_1} Y_{2}^{\beta_2} Y_{3}^{\beta_3} X_0 $\\
	$\sim (2 \alpha_1 + 2 \alpha_2 + 2 \alpha_3 - \beta_1 - \beta_2 - \beta_3) n  X_{1}^{\alpha_1} X_{2}^{\alpha_2} X_{3}^{\alpha_3}  Y_{1}^{\beta_1} Y_{2}^{\beta_2} Y_{3}^{\beta_3} X_0 $\\
	$= (2 X_1 \frac{\partial }{\partial X_1} + 2 X_2 \frac{\partial }{\partial X_2} + 2 X_3 \frac{\partial }{\partial X_3} - Y_1 \frac{\partial }{\partial Y_1} - Y_2 \frac{\partial }{\partial Y_2} - Y_3 \frac{\partial }{\partial Y_3} ) \tau_{1}^{(3)} $.\\
	
	Now let us as usual suppose $i > j$ and then we will define the following  quantities. \\
	Here for $X_i$s we have:\\
	\hspace*{3.5cm} $ _{X_j}D_{X_i} := \{ X_i , X_{j}^{n} \} $\\
	\hspace*{4.8cm}$ = X_i X_{j}^{n} - q^{2n} X_{j}^{n} X_i $\\
	\hspace*{4.8cm} $= (1 - q^{4n}) X_i X_{j}^{n} $\\
	\hspace*{4.8cm} $= (1 - e^{-4n \mathfrak{h}}) X_i X_{j}^{n} $\\
	\hspace*{4.8cm}$\sim (1 - (1 - 4n \mathfrak{h})) X_i X_{j}^{n} $\\
	\hspace*{4.8cm}$=  4n \mathfrak{h} X_i X_{j}^{n} $\\
	\hspace*{4.8cm}$\sim 4n X_i X_{j}^{n} $\\
	\hspace*{4.8cm} $= 4n X_i X_j \frac{\partial X_{j}^{n} }{\partial X_j }$.\\
	And the same will be for $Y_i$s.\\
	And for the different quantities $X_i$ and $Y_j$s we have:\\
	\hspace*{1cm} First case: for $i > j$ we have\\
	\hspace*{4cm} $ _{Y_j}D_{X_i} := \{ X_i , Y_{j}^{n} \} $\\
	\hspace*{5.20cm}$= X_i Y_{j}^{n} - q^{-n} Y_{j}^{n} X_i $\\
	\hspace*{5.20cm}$= (1 - q^{-2n}) X_i Y_{j}^{n} $\\
	\hspace*{5.20cm}$= (1 - e^{-2n \mathfrak{h}}) X_i Y_{j}^{n} $\\
	\hspace*{5.20cm}$\sim (1 - ( 1 -2n \mathfrak{h} ))  X_i Y_{j}^{n} $\\
	\hspace*{5.20cm}$= 2n \mathfrak{h} X_i Y_{j}^{n} $\\
	\hspace*{5.20cm}$\sim 2n X_i Y_{j}^{n} $\\
	\hspace*{5.18cm} $= 2 X_i Y_j \frac{\partial Y_{j}^{n} }{\partial Y_j}   $.\\
	\hspace*{1cm} Second case: for $i \leq j$ we have\\
	According to what has just mentioned we have\\
	
	\hspace*{3.8cm} $ _{Y_j}D_{1}^{Y}  := _{Y_j}D_{Y_1} $\\
	\hspace*{5.30cm} $= 4 Y_1 Y_j \frac{\partial Y_{j}^{n}}{\partial Y_j}$.\\
	And \\
	\hspace*{4.3cm} $ _{Y_1}D_{1}^{Y}   := _{Y_1}D_{Y_1} $\\
	\hspace*{5.30cm} $= 2 Y_{1}^{2} \frac{\partial Y_{1}^{n}}{\partial Y_1}. $\\
	And in a same way we can find the desired results for $_{Y_j}D_{2}^{Y} $ and $_{Y_j}D_{3}^{Y} $. \\
	So let us define
	\begin{equation}\label{Equ22}
		\begin{cases}
			_YD_{1}^{Y} := _{Y_1}D_{1}^{Y} + ^{j<1}   _{Y_j}D_{1}^{Y}  + ^{j>1}   _{Y_j}D_{1}^{Y} ; \\
			_YD_{2}^{Y} := _{Y_2}D_{2}^{Y} + ^{j<2}   _{Y_j}D_{2}^{Y}  + ^{j>2}   _{Y_j}D_{2}^{Y} ; \\
			_YD_{3}^{Y} := _{Y_3}D_{3}^{Y} + ^{j<3}   _{Y_j}D_{3}^{Y}  + ^{j>3}   _{Y_j}D_{3}^{Y} ;
		\end{cases}
	\end{equation}
	
	And then we will have
	$$_YD_{1}^{Y}  = Y_{1}^{2} \frac{\partial }{\partial Y_1} + \sum\limits_{j<1}^{} 2 Y_1 Y_j \frac{\partial }{\partial Y_j} + 0 $$
	And\\
	$$_YD_{2}^{Y}  = Y_{2}^{2} \frac{\partial }{\partial Y_2} + \sum\limits_{j<2}^{} 2 Y_2 Y_j \frac{\partial }{\partial Y_j} + 0 $$ 
	And\\ 
	$$_YD_{3}^{Y}  = Y_{3}^{2} \frac{\partial }{\partial Y_3} + \sum\limits_{j<2}^{} 2 Y_3 Y_j \frac{\partial }{\partial Y_j} + 0 $$ 
	And finally we get\\
	\hspace*{3.45cm}$_YD_{Y}^{(3)} := _YD_1  + _YD_2 + _YD_3 $\\
	\hspace*{4.5cm}$= Y_1 (Y_1 + 2 Y_2 + 2 Y_3) \frac{\partial }{\partial Y_1} + Y_2(Y_2 + 2 Y_3) \frac{\partial }{\partial Y_2} + Y_{3}^{2} \frac{\partial }{\partial Y_3} $.\\
	For $j\geq 1$ we have\\
	\hspace*{3.4cm}$_{X_1}D_{Y_j} := \{ Y_j , X_{1}^{n} \} $\\
	\hspace*{4.6cm}$= Y_j X_{1}^{n} - q^{-n} X_{1}^{n} Y_j $\\
	\hspace*{4.6cm}$= (1 - q^{-2n}) Y_j X_{1}^{n} $\\
	\hspace*{4.6cm}$= (1 - e^{2n \mathfrak{h}})  Y_j X_{1}^{n} $\\
	\hspace*{4.6cm}$\sim (1 -(1 + 2n \mathfrak{h})) Y_j X_{1}^{n} $\\
	\hspace*{4.6cm}$= - 2 n \mathfrak{h} Y_j X_{1}^{n} $\\
	\hspace*{4.6cm}$ \sim -2n Y_j X_{1}^{n} $\\
	\hspace*{4.6cm}$= -2 Y_j X_{1} \frac{\partial}{\partial X_1} $\\
	\hspace*{4.6cm}$= -2 X_1 (Y_1 + Y_2 + Y_3) \frac{\partial}{\partial X_1}. $\\
	For $j\geq 2$ we have\\
	\hspace*{3.4cm}$_{X_2}D_{Y_j} := \{ Y_j , X_{2}^{n} \} $\\
	\hspace*{4.6cm}$= Y_j X_{2}^{n} - q^{-n} X_{2}^{n} Y_j $\\
	\hspace*{4.6cm}$= (1 - q^{-2n}) Y_j X_{2}^{n} $\\
	\hspace*{4.6cm} $= (1 - e^{2n \mathfrak{h}})  Y_j X_{2}^{n} $\\
	\hspace*{4.6cm}$\sim (1 -(1 + 2n \mathfrak{h})) Y_j X_{2}^{n} $\\
	\hspace*{4.6cm}$= - 2 n \mathfrak{h} Y_j X_{2}^{n} $\\
	\hspace*{4.6cm}$\sim -2n Y_j X_{2}^{n} $\\
	\hspace*{4.6cm}$= -2 Y_j X_{2} \frac{\partial}{\partial X_2} $\\
	\hspace*{4.6cm}$ = -2 X_2 (Y_2 + Y_3) \frac{\partial}{\partial X_2} .  $\\
	For $j\geq 3$ we have\\
	\hspace*{3.4cm} $_{X_3}D_{Y_j} := \{ Y_j , X_{3}^{n} \} $\\
	\hspace*{4.7cm}$= Y_j X_{3}^{n} - q^{-n} X_{3}^{n} Y_j $\\
	\hspace*{4.7cm}$= (1 - q^{-2n}) Y_j X_{3}^{n} $\\
	\hspace*{4.7cm}$= (1 - e^{2n \mathfrak{h}})  Y_j X_{3}^{n} $\\
	\hspace*{4.7cm}$ \sim (1 -(1 + 2n \mathfrak{h})) Y_j X_{3}^{n} $\\
	\hspace*{4.7cm}$= - 2 n \mathfrak{h} Y_j X_{3}^{n} $\\
	\hspace*{4.7cm}$\sim -2n Y_j X_{3}^{n} $\\
	\hspace*{4.7cm}$= -2 Y_j X_{3} \frac{\partial}{\partial X_3} $\\
	\hspace*{4.7cm}$= - 2 X_{3} Y_3 \frac{\partial}{\partial X_3} .  $\\
	And after all these, let us define \\
	\hspace*{3.4cm} $_XD_{Y}^{(3)} := _{X_1}D_{Y_j}  + _{X_2}D_{Y_j} + _{X_3}D_{Y_j} $\\
	\hspace*{4.6cm} $= -2 X_1 (Y_1 + Y_2 + Y_3) \frac{\partial}{\partial X_1}  -2 X_2 (Y_2 + Y_3) \frac{\partial}{\partial X_2} $\\
	\hspace*{4.94cm} $- 2 X_{3} Y_3 \frac{\partial}{\partial X_3} $.\\
	And finally let us define \\
	\hspace*{3.8cm}  $D_{Y}^{(3)} := _YD_{Y}^{(3)} + _XD_{Y}^{(3)} $\\
	\hspace*{4.77cm} $= Y_1 (Y_1 + 2 Y_2 + 2 Y_3) \frac{\partial }{\partial Y_1} + Y_2(Y_2 + 2 Y_3) \frac{\partial }{\partial Y_2}$\\
	\hspace*{5.06cm}$ +  Y_{3}^{2} \frac{\partial }{\partial Y_3} -2 X_1 (Y_1 + Y_2 + Y_3) \frac{\partial}{\partial X_1}  -2 X_2 (Y_2  $\\
	\hspace*{5cm} $+ Y_3) \frac{\partial}{\partial X_2} - 2 X_{3} Y_3 \frac{\partial}{\partial X_3}$.\\
	\textbf{Next step:}\\
	Now let us try to find $D_{X}^{(3)}$:\\
	For $i > 1$, let us define $_{Y_1}D_{X_i} $ as  follows: \\
	\hspace*{3.8cm}$_{Y_1}D_{X_i} := \{ X_i , Y_{1}^{n} \} $\\
	\hspace*{4.8cm} $= X_i Y_{1}^{n} - q^{-n} Y_{1}^{n} X_i $\\
	\hspace*{4.8cm} $= (1 - q^{-2n}) X_i Y_{1}^{n} $\\
	\hspace*{4.8cm} $= (1 - e^{2n \mathfrak{h}}) X_i Y_{1}^{n} $\\
	\hspace*{4.8cm} $\sim (1 - (1 + 2n \mathfrak{h} )) X_i Y_{1}^{n} $\\
	\hspace*{4.8cm} $ = -2n \mathfrak{h} X_i Y_{1}^{n} $\\
	\hspace*{4.84cm}$\sim -2 n X_i Y_{1}^{n} = -2 X_i Y_{1} \frac{\partial }{\partial Y_1} $\\
	\hspace*{4.8cm} $= -2 Y_1 (X_2 + X_3) \frac{\partial }{\partial Y_1}. $\\
	For  $i > 2$ we have\\
	\hspace*{3.6cm} $_{Y_2}D_{X_i} := \{ X_i , Y_{2}^{n} \} $\\
	\hspace*{4.8cm} $= X_i Y_{2}^{n} - q^{-n} Y_{2}^{n} X_i $\\
	\hspace*{4.8cm} $ = (1 - q^{-2n}) X_i Y_{2}^{n} $\\
	\hspace*{4.8cm}$= (1 - e^{2n \mathfrak{h}}) X_i Y_{2}^{n} $\\
	\hspace*{4.8cm} $\sim (1 - (1 + 2n \mathfrak{h} )) X_i Y_{2}^{n} $\\
	\hspace*{4.8cm} $= -2n \mathfrak{h} X_i Y_{2}^{n} $\\
	\hspace*{4.8cm} $\sim -2 n X_i Y_{2}^{n} $\\
	\hspace*{4.8cm} $ = -2 X_i Y_{2} \frac{\partial }{\partial Y_2} $\\
	\hspace*{4.8cm} $= -2 Y_2  X_3 \frac{\partial }{\partial Y_2}. $\\
	For $i > 3$ we have 0.\\
	Let us again have the following definitions \\
	$$ _{Y_1}D_{2}^{X}  := _{Y_1}D_{X_2} = -2 Y_1 X_2 \frac{\partial Y_{1}^{n}}{\partial Y_1};$$\\
	$$ _{Y_1}D_{3}^{X}  := _{Y_1}D_{X_3} = -2 Y_1 X_3 \frac{\partial Y_{1}^{n}}{\partial Y_1};$$\\
	$$ _{Y_2}D_{3}^{X}  := _{Y_2}D_{X_3} = -2 Y_2 X_3 \frac{\partial Y_{1}^{n}}{\partial Y_2};$$\\
	Now let us define\\
	$$_XD_{Y}^{(3)} :=   _{Y_1}D_{2}^{X} + _{Y_1}D_{3}^{X} + _{Y_2}D_{3}^{X}  = - Y_1(X_2 + X_3) \frac{\partial }{\partial Y_1} - Y_2 X_3 \frac{\partial }{\partial Y_2} ;$$\\
	And now as before we have \\
	\hspace*{3.6cm} $ _{X_j}D_{1}^{X}  := _{X_j}D_{X_1} $\\
	\hspace*{4.8cm} $= 4 X_1 X_j \frac{\partial X_{j}^{n}}{\partial X_j}$.\\
	\hspace*{3.6cm} $ _{X_1}D_{1}^{X}   := _{X_1}D_{X_1} $\\
	\hspace*{4.8cm} $= 2 X_{1}^{2} \frac{\partial X_{1}^{n}}{\partial X_1}. $\\
	And in a same way we are able to define  $_{X_j}D_{2}^{X} $ and $_{X_j}D_{3}^{X} $. So let us define
	\begin{equation}\label{Equ23}
		\begin{cases}
			_XD_{1}^{X} := _{X_1}D_{1}^{X} + ^{j<1}   _{X_j}D_{1}^{X}  + ^{j>1}   _{X_j}D_{1}^{X} ; \\
			_XD_{2}^{X} := _{X_2}D_{2}^{X} + ^{j<2}   _{X_j}D_{2}^{X}  + ^{j>2}   _{X_j}D_{2}^{X} ; \\
			_XD_{3}^{X} := _{X_3}D_{3}^{X} + ^{j<3}   _{X_j}D_{3}^{X}  + ^{j>3}   _{X_j}D_{3}^{X} ;
		\end{cases}
	\end{equation}
	Then we will have
	$$\hspace*{0.85cm} _XD_{1}^{X}  = X_{1}^{2} \frac{\partial }{\partial X_1} + \sum\limits_{j<1}^{} 2 X_1 X_j \frac{\partial }{\partial X_j} + 0 $$
	And\\
	$$\hspace*{0.85cm} _XD_{2}^{X}  = X_{2}^{2} \frac{\partial }{\partial X_2} + \sum\limits_{j<2}^{} 2 X_2 X_j \frac{\partial }{\partial X_j} + 0 $$  
	And\\ 
	$$\hspace*{0.85cm} _XD_{3}^{X}  = X_{3}^{2} \frac{\partial }{\partial X_3} + \sum\limits_{j<2}^{} 2 X_3 Y_j \frac{\partial }{\partial X_j} + 0 $$
	So we will have\\
	\hspace*{3.6cm} $_XD_{X}^{(3)} := _XD_1  + _XD_2 + _XD_3 $\\
	\hspace*{4.74cm} $ = X_1 (X_1 + 2 X_2 + 2 X_3) \frac{\partial }{\partial X_1} + X_2(X_2 + 2 X_3) \frac{\partial }{\partial X_2} $\\
	\hspace*{5.15cm} $+ X_{3}^{2} \frac{\partial }{\partial X_3} $.\\
	
	And therefore as in (\ref{Equ11}) we will have the following system of $PDE$s
	\begin{equation}\label{Equ24}
		\begin{cases}
			(X_1 (X_1 + 2 X_2 + 2 X_3) \frac{\partial \tau_{1}^{(3)} }{\partial X_1} + X_2 ( X_2 + 2 X_3) \frac{\partial \tau_{1}^{(3)} }{\partial X_2} +  X_{3}^{2} \frac{\partial \tau_{1}^{(3)} }{\partial X_2} \\-  Y_1 (X_1 + X_2 + X_3) \frac{\partial \tau_{1}^{(3)} }{\partial Y_1} - Y_2(X_2 + X_3)\frac{\partial f }{\partial Y_2} - Y_3 X_3 \frac{\partial \tau_{1}^{(3)} }{\partial Y_3} )  =0; \\
			(2 X_1 \frac{\partial \tau_{1}^{(3)} }{\partial X_1} + 2 X_2 \frac{\partial \tau_{1}^{(3)} }{\partial X_2} + 2 X_3 \frac{\partial \tau_{1}^{(3)} }{\partial X_3} - Y_1 \frac{\partial \tau_{1}^{(3)} }{\partial Y_1} - Y_2 \frac{\partial \tau_{1}^{(3)} }{\partial Y_2} - Y_3 \frac{\partial \tau_{1}^{(3)} }{\partial Y_3} )  =0 ;\\
			D_{Y}^{(3)} = (Y_1 (Y_1 + 2 Y_2 + 2 Y_3) \frac{\partial \tau_{1}^{(3)} }{\partial Y_1} + Y_2 ( Y_2 + 2 Y_3) \frac{\partial \tau_{1}^{(3)} }{\partial Y_2} +  Y_{3}^{2} \frac{\partial \tau_{1}^{(3)} }{\partial Y_2} \\-  Y_1 (X_1 + X_2 + X_3) \frac{\partial \tau_{1}^{(3)} }{\partial Y_1} - Y_2(X_2 + X_3)\frac{\partial \tau_{1}^{(3)} }{\partial Y_2} - Y_3 X_3 \frac{\partial \tau_{1}^{(3)} }{\partial Y_3} )  =0; \\
			(2 X_1 \frac{\partial \tau_{1}^{(3)} }{\partial X_1} + 2 X_2 \frac{\partial \tau_{1}^{(3)} }{\partial X_2} + 2 X_3 \frac{\partial \tau_{1}^{(3)} }{\partial X_3} - Y_1 \frac{\partial \tau_{1}^{(3)} }{\partial Y_1} - Y_2 \frac{\partial \tau_{1}^{(3)} }{\partial Y_2} - Y_3 \frac{\partial \tau_{1}^{(3)} }{\partial Y_3} )  =0 ;
		\end{cases}
	\end{equation}
	
	And according to appendix A, we have the following functional dependent nontrivial solution for the whole system of $PDE$s  (\ref{Equ24}) 
	
	\begin{equation}\label{Equ25}
		\hspace{-1.05cm}\tau_{1}^{(3)} = \frac{(\Sigma_{1 \leq i \leq j \leq 2} ~ ~  ~  ~ X_i Y_j)(\Sigma_{1 \leq i \leq j \leq 2} ~ ~  ~  ~ X_{i+1} Y_{j+1})}{X_2 Y_2 (\Sigma_{1 \leq i \leq j \leq 3} ~ ~  ~  ~ X_i Y_j)};
	\end{equation}
	And again as before, $(3)$ goes back to 3 in the $Sl_3$ and  $1$ is a default index which later we will use it for to employ our shifting operators.\\
	According to the number of variables, we will have 6 shifts and then after that it will be in a loop.\\
	So here in $sl_3$ case we have six solutions  which belong to the fraction ring of polynomial functions.
	\begin{equation}\label{Equ26}
		\begin{cases}
			
			\tau_{1}^{(3)}[X_1, Y_1, X_2, Y_2, X_3, Y_3]= \frac{X_2 Y_2 (X_3 Y_3 + X_2 (Y_2 + Y_3) + X_1 (Y_1 + Y_2 + Y_3))}{(X_2 Y_2 + X_1 (Y_1 + Y_2)) (X_3 Y_3 + X_2 (Y_2 + Y_3))} ; \\
			
			\vspace{0.2cm}
			
			\tau_{2}^{(3)}[Y_1, X_2, Y_2, X_3, Y_3, X_4]= \frac{X_3 Y_2 (X_2 Y_1 + (X_3 + X_4) (Y_1 + Y_2) + X_4 Y_3)}{(X_2 Y_1 + X_3 (Y_1 + Y_2)) (X_3 Y_2 + X_4 (Y_2 + Y_3))}  ; \\
			\vspace{0.2cm}
			
			\tau_{3}^{(3)}[X_2, Y_2, X_3, Y_3, X_4, Y_4]= \frac{X_3 y_3 (X_4 Y_4 + X_3 (Y_3 + Y_4) + X_2 (Y_2 + Y_3 + Y_4))}{ (X_3 Y_3 + X_2 (Y_2 + Y_3)) (X_4 Y_4 + X_3 (Y_3 + Y_4)) } ; \\
			\vspace{0.2cm}
			\tau_{4}^{(3)}[Y_2, X_3, Y_3, X_4, Y_4, X_5]= \frac{X_4 Y_3 (X_3 Y_2 + (X_4 + X_5) (Y_2 + Y_3) + X_5 Y_4)}{(X_3 Y_2 + X_4 (Y_2 + Y_3)) (X_4 Y_3 + X_5 (Y_3 + Y_4))} ;\\
			\vspace{0.2cm}
			\tau_{5}^{(3)}[X_3, Y_3, X_4, Y_4, X_5, Y_5]=  \frac{X_4 Y_4 (X_5 Y_5 + X_4 (Y_4 + Y_5) + X_3 (Y_3 + Y_4 + Y_5))}{(X_4 Y_4 + X_3 (Y_3 + Y_4)) (X_5 Y_5 + X_4 (Y_4 + Y_5))} ;\\
			\vspace{0.2cm}
			\tau_{6}^{(3)}[Y_3, X_4, Y_4, X_5, Y_5, X_6]= \frac{X_5 Y_4 (X_4 Y_3 + (X_5 + X_6) (Y_3 + Y_4) + X_6 Y_5)}{(X_4 Y_3 + X_5 (Y_3 + Y_4)) (X_5 Y_4 + X_6 (Y_4 + Y_5))};
		\end{cases}
	\end{equation}
	Where $\tau_{1}^{(3)} := \tau_{1}^{(3)}[\cdots  , X_1, Y_1 , X_2, Y_2,  X_3, Y_3  \cdots ]$. \\
	Again by setting $X_{i}^{(1i)} := X_i$ and $X_{i}^{(2i)} := Y_i$ and $X_{i}^{(3i)} := Z_i$ and according to $(\ref{Equ17})$ we have to write down the following brackets as a composition of $\tau_{i}^{(3)}$s, because of the algebra structure and it will be done by using Mathematica coding in appendix A. 
	\begin{equation}\label{Equ27}
		\begin{cases}
			
			F_{2}^{(3)} = \{ \tau_{1}^{(3)} , \tau_{2}^{(3)} \} = -(1 - \tau_{1}^{(3)})(1 - \tau_{2}^{(3)}) (\tau_{1}^{(3)} \tau_{2}^{(3)}); \\
			
			\vspace{0.2cm}
			
			F_{3}^{(3)} = \{ \tau_{1}^{(3)} , \tau_{3}^{(3)} \} = (1 - \tau_{1}^{(3)})(1 - \tau_{3}^{(3)}) (\tau_{1}^{(3)} \tau_{2}^{(3)} + \tau_{2}^{(3)} \tau_{3}^{(3)} - \tau_{2}^{(3)} ) ; \\
			\vspace{0.2cm}
			
			F_{4}^{(3)} = \{ \tau_{1}^{(3)} , \tau_{4}^{(3)} \} = -(1 - \tau_{1}^{(3)})(1 - \tau_{4}^{(3)})\\ \hspace*{1.3cm} (\tau_{1}^{(3)}  \tau_{2}^{(3)} + \tau_{2}^{(3)} \tau_{3}^{(3)} + \tau_{3}^{(3)} \tau_{4}^{(3)} - \tau_{1}^{(3)} - \tau_{2}^{(3)} - \tau_{3}^{(3)} - \tau_{4}^{(3)} + 1) ; \\
			\vspace{0.2cm}
			F_{5}^{(3)} = \{ \tau_{1}^{(3)} , \tau_{5}^{(3)} \} = (1 - \tau_{1}^{(3)})(1 - \tau_{5}^{(3)}) (\tau_{2}^{(3)} \tau_{3}^{(3)} + \tau_{3}^{(3)} \tau_{4}^{(3)} - \tau_{2}^{(3)} \\ \hspace*{1.3cm} - \tau_{3}^{(3)} - \tau_{4}^{(3)} + 1 );\\
			\vspace{0.2cm}
			F_{6}^{(3)} = \{ \tau_{1}^{(3)} , \tau_{6}^{(3)} \} = -(1 - \tau_{1}^{(3)})(1 - \tau_{6}^{(3)})( \tau_{3}^{(3)} \tau_{4}^{(3)}  - \tau_{4}^{(3)} - \tau_{3}^{(3)} + 1 ) ;\\
			\vspace{0.2cm}
			F_{i}^{(3)} = \{ \tau_{1}^{(3)} , \tau_{i}^{(3)} \} = 0 ~~~~~ \hspace*{3cm} ~~~ \ ~~~~  \text{for } |i - 1| \geq 6 ;\\
		\end{cases}
	\end{equation}
	\subsection{Lattice $W_4$ algebra; main generator}
	In this case we will use the following defined Poisson bracket (\ref{Equ28}) based on Cartan matrix \\
	\hspace*{4cm}$A_3 =  \left[ \begin{array}{ccc}
	$2$&$-1$&$0$ \\
	$-1$&$2$&$-1$ \\
	$0$&$-1$&$2$ \\
	\end{array} \right].$ \\
	But for to do this according to our previous ordering and list of variables, and the same as what we have done in $sl_3$ case, let us for simplicity of the calculations, order the set of our variables as follows:\\
	Set $X_{i}^{(1i)} := X_i$ and $X_{i}^{(2i)} := Y_i$ and $X_{i}^{(3i)} := Z_i$ and so on.\\
	\begin{defn}
		Let's define our Poisson bracket as follows in the case of $sl_4$:
		\begin{equation}\label{Equ28}
			\begin{cases}
				\{X_i , X_j \} := 2 X_i X_j & \text{if } i < j;\\
				\{Y_i , Y_j \} := 2 Y_i Y_j & \text{if } i < j;\\
				\{Z_i , Z_j \} := 2 Z_i Z_j & \text{if } i < j;\\
				\{ X_i , X_i \} := 0; \\
				\{ Y_i , Y_i \} := 0; \\
				\{ Z_i , Z_i \} := 0; \\
				\{ X_i , Y_j \} := X_i Y_j & \text{if } i > j;\\
				\{ X_i , Y_j \} := - X_i Y_j & \text{if } i \leq j;\\
				\{ X_i , Z_j \} := 0; \\
				\{ Y_i , Z_j \} := Y_i Z_j & \text{if } i > j;\\
				\{ Y_i , Z_j \} := - Y_i Z_j & \text{if } i \leq j;\\
			\end{cases}
		\end{equation}

	\end{defn}
	
	
	And instead of (\ref{Equ1}) we will have the following $q-$commutation relations for $j \in \{ 1, 2, 3 \}$ and as always $i \in \{ 1, 2 , 3 \}$:
	\begin{equation}\label{Equ29}
		\begin{cases}
			X_i X_j = q^2 X_j X_i& \text{if } i \leq j\\
			Y_i Y_j = q^2 Y_j Y_i& \text{if } i \leq j \\
			Z_i Z_j = q^2 Z_j Z_i& \text{if } i \leq j \\
			X_i Y_j = q^{-1} Y_j X_i& \text{if } i \leq j \\
			Y_i Z_j = q^{-1} Z_j Y_i& \text{if } i \leq j \\
			X_i Z_j  = Z_j X_i& \text{ } 
			
		\end{cases}
	\end{equation}

	And by using the same approach as in the  $sl_2$ and $sl_3$ case, it become clear that  the equations $D_{X}^{(4)}$, $D_{Y}^{(4)}$ and $D_{Z}^{(4)}$ and also $H_{X}^{(4)}$, $H_{Y}^{(4)}$ and $H_{Z}^{(4)}$    will have the following forms:
	\begin{equation}\label{Equ30}
		{\mathfrak{D}}_{X}^{(4)} = X_1 (X_1 + 2 X_2 + 2 X_3) \frac{\partial \tau_{1}^{(4)}}{\partial X_1} + X_2 (X_2 + 2 X_3) \frac{\partial \tau_{1}^{(4)}}{\partial X_2} + X_{3}^{2} \frac{\partial \tau_{1}^{(4)}}{\partial X_3} - Y_1 (X_2 + X_3)
	\end{equation}
	$$\hspace{-7.15cm}  \frac{\partial \tau_{1}^{(4)}}{\partial Y_1}  - Y_2 X_3 \frac{\partial \tau_{1}^{(4)}}{\partial Y_2};$$
	\begin{equation}\label{Equ31}
		\hspace{-0.06cm}{\mathfrak{D}}_{Y}^{(4)} = Y_1 (Y_1 + 2 Y_2 + 2 Y_3) \frac{\partial \tau_{1}^{(4)}}{\partial Y_1} + Y_2 (Y_2 + 2 Y_3) \frac{\partial \tau_{1}^{(4)}}{\partial Y_2} + Y_{3}^{2} \frac{\partial \tau_{1}^{(4)}}{\partial Y_3} - X_1 (Y_1 + Y_2 + Y_3)  
	\end{equation}
	$$\hspace{1.1cm} \frac{\partial \tau_{1}^{(4)}}{\partial X_1} - X_2( Y_2 + Y_3 ) \frac{\partial \tau_{1}^{(4)}}{\partial X_2} - X_3 Y_3 \frac{\partial \tau_{1}^{(4)}}{\partial X_3} - Z_1 (Y_2 + Y_3) \frac{\partial \tau_{1}^{(4)}}{\partial z_1} - Z_2 y_3 \frac{\partial \tau_{1}^{(4)}}{\partial Z_2};$$
	\begin{equation}\label{Equ32}
		\hspace{-0.07cm}{\mathfrak{D}}_{Z}^{(4)} =Z_1 (Z_1 + 2 Z_2 + 2 Z_3) \frac{\partial \tau_{1}^{(4)}}{\partial Z_1} + Z_2 (Z_2 + 2 Z_3) \frac{\partial \tau_{1}^{(4)}}{\partial Z_2} + Z_{3}^{2} \frac{\partial \tau_{1}^{(4)}}{\partial Z_3} - Y_1 (Z_1 + Z_2 + Z_3) 
	\end{equation}
	$$\hspace{-3.6cm} \frac{\partial \tau_{1}^{(4)}}{\partial Y_1}  - Y_2 (Z_2 + Z_3) \frac{\partial \tau_{1}^{(4)}}{\partial Y_2} - Y_3 Z_3 \frac{\partial \tau_{1}^{(4)}}{\partial Y_3};$$
	\begin{equation}\label{Equ33}
		\hspace{-0.1cm}{H}_{X}^{(4)} = 2 X_1 \frac{\partial \tau_{1}^{(4)}}{\partial X_1} +  2 X_2 \frac{\partial \tau_{1}^{(4)}}{\partial X_2} +  2 X_3 \frac{\partial \tau_{1}^{(4)}}{\partial X_3} -  Y_1 \frac{\partial \tau_{1}^{(4)}}{\partial Y_1} -  Y_2 \frac{\partial \tau_{1}^{(4)}}{\partial Y_2} -  Y_3 \frac{\partial \tau_{1}^{(4)}}{\partial Y_3};
	\end{equation}
	\begin{equation}\label{Equ34}
		\hspace{-0.14cm}{H}_{Y}^{(4)} = 2 Y_1 \frac{\partial \tau_{1}^{(4)}}{\partial Y_1} +  2 Y_2 \frac{\partial \tau_{1}^{(4)}}{\partial Y_2} +  2 Y_3 \frac{\partial \tau_{1}^{(4)}}{\partial Y_3} -  X_1 \frac{\partial \tau_{1}^{(4)}}{\partial X_1} -  X_2 \frac{\partial \tau_{1}^{(4)}}{\partial X_2} -  X_3 \frac{\partial \tau_{1}^{(4)}}{\partial X_3} -  z_1 \frac{\partial \tau_{1}^{(4)}}{\partial Z_1}  
	\end{equation}
	$$\hspace{-7.0cm} -  Z_2 \frac{\partial \tau_{1}^{(4)}}{\partial Z_2} -Z_3 \frac{\partial \tau_{1}^{(4)}}{\partial Z_3};$$
	\begin{equation}\label{Equ35}
		\hspace{-0.1cm}{H}_{Z}^{(4)} = 2 Z_1 \frac{\partial \tau_{1}^{(4)}}{\partial Z_1} +  2 Z_2 \frac{\partial \tau_{1}^{(4)}}{\partial Z_2} +  2 Z_3 \frac{\partial \tau_{1}^{(4)}}{\partial Z_3} -  Y_1 \frac{\partial \tau_{1}^{(4)}}{\partial Y_1} -  Y_2 \frac{\partial \tau_{1}^{(4)}}{\partial Y_2} -  Y_3 \frac{\partial \tau_{1}^{(4)}}{\partial Y_3};
	\end{equation}
	And the functional dependent nontrivial solutions for the whole system of first order partial differential equation is as follows:
	
	\begin{equation}\label{Equ36}
		\tau_{1}^{(4)} = \frac{(\Sigma_{1 \leq i \leq j \leq m \leq 2} ~ ~  ~  ~ x_i y_j z_m)(\Sigma_{1 \leq i \leq j \leq m \leq 2} ~ ~  ~  ~ x_{i+1} y_{j+1} z_{m+1})}{x_2 y_2 z_2 (\Sigma_{1 \leq i \leq j \leq m \leq 3} ~ ~  ~  ~ x_i y_j z_m)};
	\end{equation}
	And again as before, $(4)$ goes back to 4 in the $Sl_4$ and  $1$ is a default index which later we will use it for to employ our shifting operators.\\
	According to the number of variables, we will have 9 shifts and then after that it will be in a loop.\\
	So here in $sl_4$ case we have nine solutions: \\

	\hspace*{-0.4cm}$\tau_{1}^{(4)} := \tau_{1}^{(4)}[X_1, Y_1, Z_1, X_2, Y_2, Z_2, X_3, Y_3, Z_3]$;

	\hspace*{-0.4cm}$\tau_{2}^{(4)} := \tau_{1}^{(4)}[X_1 \rightarrow Y_1, Y_1 \rightarrow Z_1, Z_1 \rightarrow  X_2 , X_2 \rightarrow Y_2, Y_2 \rightarrow Z_2, Z_2 \rightarrow  X_3 , X_3 \rightarrow Y_3 , Y_3 \rightarrow Z_3 ];$

	\hspace*{-0.4cm}$\tau_{3}^{(4)} := \tau_{2}^{(4)}[Y_1 \rightarrow Z_1, Z_1 \rightarrow X_2, X_2 \rightarrow  Y_2 , Y_2 \rightarrow Z_2, Z_2 \rightarrow  X_3, X_3 \rightarrow Y_3 , Y_3 \rightarrow Z_3 , Z_3 \rightarrow X_4 ];$

	\hspace*{-0.4cm}$\tau_{4}^{(4)} := \tau_{3}^{(4)}[Z_1 \rightarrow X_2, X_2 \rightarrow Y_2, Y_2 \rightarrow  Z_2 , Z_2 \rightarrow   X_3, X_3 \rightarrow Y_3, Y_3 \rightarrow Z_3 , Z_3 \rightarrow X_4 , X_4 \rightarrow Y_4 ];$

	\hspace*{-0.4cm}$\tau_{5}^{(4)} := \tau_{4}^{(4)}[X_2 \rightarrow Y_2, Y_2 \rightarrow Z_2, Z_2 \rightarrow   X_3 , X_3 \rightarrow Y_3, Y_3 \rightarrow Z_3, Z_3 \rightarrow X_4 , X_4 \rightarrow  Y_4 , Y_4 \rightarrow Z_4 ];$

	\hspace*{-0.4cm}$\tau_{6}^{(4)} := \tau_{5}^{(4)}[Y_2 \rightarrow Z_2, Z_2 \rightarrow  X_3, X_3 \rightarrow  Y_3 , Y_3 \rightarrow Z_3, Z_3 \rightarrow X_4, X_4 \rightarrow Y_4 , Y_4 \rightarrow  Z_4 , Z_4 \rightarrow X_5 ];$

	\hspace*{-0.4cm}$\tau_{7}^{(4)} := \tau_{6}^{(4)}[Z_2 \rightarrow   X_3, X_3 \rightarrow Y_3, Y_3 \rightarrow  Z_3 , Z_3 \rightarrow X_4, X_4 \rightarrow Y_4, Y_4 \rightarrow Z_4 , Z_4 \rightarrow   X_5 , X_5 \rightarrow Y_5 ];$

	\hspace*{-0.4cm}$\tau_{8}^{(4)} := \tau_{7}^{(4)}[X_3 \rightarrow Y_3, Y_3 \rightarrow Z_3, Z_3 \rightarrow  X_4 , X_4 \rightarrow Y_4, Y_4 \rightarrow Z_4, Z_4 \rightarrow X_5 , X_5 \rightarrow  Y_5 , Y_5 \rightarrow Z_5 ];$

	\hspace*{-0.4cm}$\tau_{9}^{(4)} := \tau_{8}^{(4)}[Y_3 \rightarrow Z_3, Z_3 \rightarrow X_4, X_4 \rightarrow  Y_4 , Y_4 \rightarrow Z_4, Z_4 \rightarrow X_5, X_5 \rightarrow Y_5 , Y_5 \rightarrow Z_5 , Z_5 \rightarrow X_6 ];$
	
	which belong to the fraction ring of polynomial functions.
	\subsection{Lattice $W_5$ algebra; main generator}
	In this case we will use the following defined Poisson bracket based on Cartan matrix \\
	\hspace*{4cm}$A_4 =  \left[ \begin{array}{cccc}
	$2$&$-1$&$0$&$0$ \\
	$-1$&$2$&$-1$&$0$ \\
	$0$&$0$&$-1$&$2$ \\
	\end{array} \right].$ \\
	But for to do this according to our previous ordering and list of variables, and the same as what we have done in $sl_4$ case, let us for simplicity in the calculations,  order the set of our variables as follows:\\
	Set $X_{i}^{(1i)} := X_i$ and $X_{i}^{(2i)} := Y_i$ and $X_{i}^{(3i)} := Z_i$ and $X_{i}^{(4i)} := K_i$.\\
	\begin{defn}
		Let's define our Poisson bracket as follows in the case of $sl_5$:
		\begin{equation}\label{Equ37}
			\begin{cases}
				\{X_i , X_j \} := 2 X_i X_j & \text{if } i < j;\\
				\{Y_i , Y_j \} := 2 Y_i Y_j & \text{if } i < j;\\
				\{Z_i , Z_j \} := 2 Z_i Z_j & \text{if } i < j;\\
				\{K_i , K_j \} := 2 K_i K_j & \text{if } i < j;\\
				\{ X_i , X_i \} := 0; \\
				\{ Y_i , Y_i \} := 0; \\
				\{ Z_i , Z_i \} := 0; \\
				\{ K_i , K_i \} := 0; \\
				\{ X_i , Y_j \} := X_i Y_j & \text{if } i > j;\\
				\{ X_i , Y_j \} := - X_i Y_j & \text{if } i \leq j;\\
				\{ X_i , Z_j \} := 0; \\
				\{ X_i , K_j \} := 0; \\
				\{ Y_i , Z_j \} := Y_i Z_j & \text{if } i > j;\\
				\{ Y_i , Z_j \} := - Y_i Z_j & \text{if } i \leq j;\\
				\{ Y_i , K_j \} := Y_i K_j & \text{if } i > j;\\
				\{ Y_i , K_j \} := - Y_i K_j & \text{if } i \leq j;\\
			\end{cases}
		\end{equation}

	\end{defn}
	
	
	And instead of $(\ref{Equ1})$ we will have the following $q-$commutation relations for $j \in \{ 1, 2, 3 \}$ and as always $i \in \{ 1, 2 , 3 \}$:
	\[ 
	\hspace{-1cm}\begin{cases}
	x_i x_j = q^2 x_j x_i& \text{if } i \leq j\\
	y_i y_j = q^2 y_j y_i& \text{if } i \leq j \\
	z_i z_j = q^2 z_j z_i& \text{if } i \leq j \\
	k_i k_j = q^2 k_j k_i& \text{if } i \leq j \\
	x_i y_j = q^{-1} y_j x_i& \text{if } i \leq j \\
	y_i z_j = q^{-1} z_j y_i& \text{if } i \leq j \\
	z_i k_j = q^{-1} k_j z_i& \text{if } i \leq j \\
	x_i z_j  = z_j x_i& \text{ } \\
	y_i k_j = k_j y_i& \text{ } \\
	x_i k_j = k_j x_i& \text{ }
	\end{cases}
	\]
	
	And by using the same approach as  $sl_2$, $sl_3$ and $sl_4$ case , it become clear that  the equations $D_{X}^{(5)}$, $D_{Y}^{(5)}$, $D_{Z}^{(5)}$ and $D_{K}^{(5)}$ and also $H_{X}^{(5)}$, $H_{Y}^{(5)}$, $H_{Z}^{(5)}$ and $H_{K}^{(5)}$     will have the following forms:
	\begin{equation}\label{Equ38}
		\hspace{-0.05cm}{\mathfrak{D}}_{X}^{(5)} = X_1 (X_1 + 2 X_2 + 2 X_3) \frac{\partial \tau_{1}^{(5)}}{\partial X_1} + X_2 (X_2 + 2 X_3) \frac{\partial \tau_{1}^{(5)}}{\partial X_2} + X_{3}^{2} \frac{\partial \tau_{1}^{(5)}}{\partial X_3} - Y_1 (X_2 + X_3)  
	\end{equation}
	$$ \hspace{-7.35cm} \frac{\partial \tau_{1}^{(5)}}{\partial Y_1} - Y_2 X_3 \frac{\partial \tau_{1}^{(5)}}{\partial Y_2};$$
	\begin{equation}\label{Equ39}
		\hspace{-0.05cm}{\mathfrak{D}}_{Y}^{(5)} = Y_1 (Y_1 + 2 Y_2 + 2 Y_3) \frac{\partial \tau_{1}^{(5)}}{\partial Y_1} + Y_2 (Y_2 + 2 Y_3) \frac{\partial \tau_{1}^{(5)}}{\partial Y_2} + Y_{3}^{2} \frac{\partial \tau_{1}^{(5)}}{\partial Y_3} - X_1 (Y_1 + Y_2 + Y_3) 
	\end{equation}
	$$\hspace{0.85cm} \frac{\partial \tau_{1}^{(5)}}{\partial X_1}  - X_2( Y_2 + Y_3 ) \frac{\partial \tau_{1}^{(5)}}{\partial X_2} - X_3 Y_3 \frac{\partial \tau_{1}^{(5)}}{\partial X_3} - Z_1 (Y_2 + Y_3) \frac{\partial \tau_{1}^{(5)}}{\partial z_1} - Z_2 y_3 \frac{\partial \tau_{1}^{(5)}}{\partial Z_2};$$
	\begin{equation}\label{Equ40}
		\hspace{-0.05cm}{\mathfrak{D}}_{Z}^{(5)} =Z_1 (Z_1 + 2 Z_2 + 2 Z_3) \frac{\partial \tau_{1}^{(5)}}{\partial Z_1} + Z_2 (Z_2 + 2 Z_3) \frac{\partial \tau_{1}^{(5)}}{\partial Z_2} + Z_{3}^{2} \frac{\partial \tau_{1}^{(5)}}{\partial Z_3} - Y_1 (Z_1 + Z_2 + Z_3)  
	\end{equation}
	$$\hspace{0.85cm} \frac{\partial \tau_{1}^{(5)}}{\partial Y_1} - Y_2 (Z_2 + Z_3) \frac{\partial \tau_{1}^{(5)}}{\partial Y_2} - Y_3 Z_3 \frac{\partial \tau_{1}^{(5)}}{\partial Y_3}  - K_1 (Z_2 + Z_3  ) \frac{\partial \tau_{1}^{(5)}}{\partial k_1} - 
	K_2 Z_3  \frac{\partial \tau_{1}^{(5)}}{\partial K_2};$$
	\begin{equation}\label{Equ41}
		\hspace{-0.15cm}{\mathfrak{D}}_{K}^{(5)} =K_1 (K_1 + 2 K_2 + 2 K_3) \frac{\partial \tau_{1}^{(5)}}{\partial K_1} + K_2 (K_2 + 2 K_3) \frac{\partial \tau_{1}^{(5)}}{\partial K_2} + K_{3}^{2} \frac{\partial \tau_{1}^{(5)}}{\partial Z_3} - Z_1 (K_1 + K_2 
	\end{equation}
	$$\hspace{-3.06cm} + K_3) \frac{\partial \tau_{1}^{(5)}}{\partial Z_1}  - Z_2 ( K_2 + K_3) \frac{\partial \tau_{1}^{(5)}}{\partial z_2}  - Z_3 X_3 \frac{\partial \tau_{1}^{(5)}}{\partial Z_3};$$
	\begin{equation}\label{Equ42}
		{H}_{X}^{(5)} = 2 X_1 \frac{\partial \tau_{1}^{(5)}}{\partial X_1} +  2 X_2 \frac{\partial \tau_{1}^{(5)}}{\partial X_2} +  2 X_3 \frac{\partial \tau_{1}^{(5)}}{\partial X_3} -  Y_1 \frac{\partial \tau_{1}^{(5)}}{\partial Y_1} -  Y_2 \frac{\partial \tau_{1}^{(5)}}{\partial Y_2} -  Y_3 \frac{\partial \tau_{1}^{(5)}}{\partial Y_3};
	\end{equation}
	\begin{equation}\label{Equ43}
		{H}_{Y}^{(5)} = 2 Y_1 \frac{\partial \tau_{1}^{(5)}}{\partial Y_1} +  2 Y_2 \frac{\partial \tau_{1}^{(5)}}{\partial Y_2} +  2 Y_3 \frac{\partial \tau_{1}^{(5)}}{\partial Y_3} -  X_1 \frac{\partial \tau_{1}^{(5)}}{\partial X_1} -  X_2 \frac{\partial \tau_{1}^{(5)}}{\partial X_2} -  X_3 \frac{\partial \tau_{1}^{(5)}}{\partial X_3} -  Z_1 \frac{\partial \tau_{1}^{(5)}}{\partial Z_1} 
	\end{equation}
	$$\hspace{-6.36cm} - Z_2 \frac{\partial \tau_{1}^{(5)}}{\partial Z_2} -  Z_3 \frac{\partial \tau_{1}^{(5)}}{\partial Z_3};$$
	\begin{equation}\label{Equ44}
		{H}_{Z}^{(5)} = 2 Z_1 \frac{\partial \tau_{1}^{(5)}}{\partial Z_1} +  2 Z_2 \frac{\partial \tau_{1}^{(5)}}{\partial Z_2} +  2 Z_3 \frac{\partial \tau_{1}^{(5)}}{\partial z_3} -  Y_1 \frac{\partial \tau_{1}^{(5)}}{\partial Y_1} -  y_2 \frac{\partial \tau_{1}^{(5)}}{\partial Y_2} -  Y_3 \frac{\partial \tau_{1}^{(5)}}{\partial Y_3} -  K_1 \frac{\partial \tau_{1}^{(5)}}{\partial K_1} 
	\end{equation}
	$$\hspace{-6.36cm}  - K_2 \frac{\partial \tau_{1}^{(5)}}{\partial K_2} -  K_3 \frac{\partial \tau_{1}^{(5)}}{\partial K_3};$$
	\begin{equation}\label{Equ45}
		{H}_{K}^{(5)} = 2 K_1 \frac{\partial \tau_{1}^{(5)}}{\partial K_1} +  2 K_2 \frac{\partial \tau_{1}^{(5)}}{\partial K_2} +  2 K_3 \frac{\partial \tau_{1}^{(5)}}{\partial K_3} -  Z_1 \frac{\partial \tau_{1}^{(5)}}{\partial Z_1} -  Z_2 \frac{\partial \tau_{1}^{(5)}}{\partial Z_2} -  Z_3 \frac{\partial \tau_{1}^{(5)}}{\partial Z_3};
	\end{equation}
	
	And the functional dependent nontrivial solution for the whole system of first order partial differential equation is as follows:
	
	\begin{equation}\label{Equ46}
		\tau_{1}^{(5)} = \frac{(\Sigma_{1 \leq i \leq j \leq m \leq l \leq 2} ~ ~  ~  ~ x_i y_j z_m k_l )(\Sigma_{1 \leq i \leq j \leq m \leq l \leq 2} ~ ~  ~  ~ x_{i+1} y_{j+1} z_{m+1} k_{l+1})}{x_2 y_2 z_2 k_2 (\Sigma_{1 \leq i \leq j \leq m \leq l \leq 3} ~ ~  ~  ~ x_i y_j z_m k_l)};
	\end{equation}
	And again as before, $(5)$ goes back to 5 in the $Sl_5$ and  $1$ is a default index which later we will use it for to employ our shifting operators.\\
	According to the number of variables, we will have 12 shifts and then after that it will be in a loop.\\
	So here in $sl_5$ case we will have twelve solutions just as what we did in the $sl_4$ case, and here we skip to write them down. \\
	\subsection{Lattice $W_n$ algebra; main generator}
	Here for $sl_n$, we skip writing down all steps which we have been done in the previous sections and we just will write down the main generator of the lattice $W_n$ algebra.\\
	The functional dependent nontrivial solution for the whole system of the first order partial differential equations will be as what comes in follow:
	\begin{equation}\label{Equ47}
		\hspace{-1.5cm} \tau_{1}^{(n)} = 
	\end{equation}
	$$\frac{(\Sigma_{1 \leq i_1 \leq i_2 \cdots \leq i_{n-1} \leq 2} ~ ~  ~  ~ x_{i_1}^{(1)} x_{i_2}^{(2)} \cdots x_{i_{n-1}}^{(n-1)} )(\Sigma_{1 \leq i_1 \leq i_2 \cdots \leq i_{n-1} \leq 2} ~ ~  ~  ~ x_{i_1 +1}^{(1)} x_{i_2 +1}^{(2)} \cdots x_{i_{n-1} +1}^{(n-1)})}{x_{2}^{(1)} \cdots x_{2}^{(n-1)} (\Sigma_{1 \leq i_1 \leq i_2 \cdots \leq i_{n-1} \leq 3} ~ ~  ~  ~ x_{i_1}^{(1)} x_{i_2}^{(2)} \cdots x_{i_{n-1}}^{(n-1)} )}.$$
	We should notice that $x_{i_j}^{(j)}$s are different of each other for any $j \in \{ 1, \cdots n-1 \} $
	
	\section*{Acknowledgement}
	The research in this paper would have taken far longer to complete without the encouragement from many others. It is a delight to acknowledge those who have supported me over the last three years during the preparation of this paper. 
	
	I would like to thank Prof. Yaroslav Pugai, for his guidance and relaxed, thoughtful insight. 
	
	I thank all of the Institute for information transmission problems (Kharkevich Institute)'s staff for their hospitality over the last year.
	
	I am particularly thankful for the help and advice of Prof. Brendan Godfrey, without whom the learning Mathematica would have been very much steeper and unimaginable.
	
	And I would like to thank Prof. Boris Feigin for suggesting me this interesting problem and enlightening discussions during the preparation for my first paper which was my first step in this subject!
	
	And finally, I would like to thank Prof. Alexei Kanel-Belov for his support and for whom I will always be indebted for being a constant source of inspiration and for the great talent and patience which he has always guided me with.

		\section{Appendix A}
		\hspace*{2.4cm} $ \textit{Brendan B. Godfrey  } $ ~~~~~ \and ~~~~~ $\textit{Farrokh Razavinia }$\\
		 \hspace*{4cm} \href{brendan.godfrey@ieee.org}{brendan.godfrey@ieee.org}   \\
		 \hspace*{4cm} \href{f.razavinia@phystech.edu}{f.razavinia@phystech.edu}  \\

		This section has been completed by getting help from professor Brendan B. Godfrey from the Institute for Research in Electronics and Applied Physics (The University of Maryland), indirect communications and discussions through email, and also through a series of questions and discussions in Mathematica StackExchange.
		
		And I have to say that without his great Mathematica skills, it nearly was impossible to get such interesting results!
		
		In this appendix, you will find some parts of Mathematica codings which we have used to obtain our algebra structures. 
		
		And we believe that what is written in this appendix can open a new approach in solving the following system of $q-$linear homogeneous equations in one unknown $f$.
		\begin{equation}
		\begin{cases}
		equ_1(f) = a_{11} \frac{\partial f}{\partial x_1} + a_{21} \frac{\partial f}{\partial x_2} + \cdots + a_{n1} \frac{\partial f}{\partial x_n} = 0 \\
		equ_2(f) = a_{12} \frac{\partial f}{\partial x_1} + a_{22} \frac{\partial f}{\partial x_2} + \cdots + a_{n2} \frac{\partial f}{\partial x_n} = 0 \\
		\hspace*{0.3 cm}  \vdots  \hspace*{1.9 cm}  \vdots \hspace*{1.3 cm} \vdots \hspace*{2.3 cm} \vdots  \\
		equ_q(f) = a_{1q} \frac{\partial f}{\partial x_1} + a_{2q} \frac{\partial f}{\partial x_2} + \cdots + a_{nq} \frac{\partial f}{\partial x_n} = 0
		\end{cases}
		\end{equation}
	
		Where the coefficients $a_{ik}$ are functions of $n$ independent variables $x_1, \cdots , x_n $ and do not contain the unknown function $f$. \cite{3}
		
		And we have to mention that, to reach to this point was impossible without using Mathematica!
	\subsection{Lattice $W_3$ algebra}.
	
	\begin{lstlisting}[language=Mathematica,caption={Example code}]
	p = D[f[x1, x2, x3, y1, y2, y3], x1];
	q = D[f[x1, x2, x3, y1, y2, y3], x2];
	r = D[f[x1, x2, x3, y1, y2, y3], x3];
	o = D[f[x1, x2, x3, y1, y2, y3], y1];
	x = D[f[x1, x2, x3, y1, y2, y3], y2];
	a = D[f[x1, x2, x3, y1, y2, y3], y3];
	equ1 = 2 x1 p + 2 x2 q + 2 x3 r - y1 o - y2 x - y3 a;
	equ2 = -x1 p - x2 q - x3 r + 2 y1 o + 2 y2 x + 2 y3 a;
	equ3 = (x1 (x1 + 2 x2 + 2 x3)) p + (x2 (x2 + 2 x3)) q + (x3^2) r 
	- (y1 (x2 + x3)) o - y2 x3 x;
	equ4 = (y1 (y1 + 2 y2 + 2 y3)) o + (y2 (y2 + 2 y3)) x + (y3^2) a 
	- (x1 (y1 + y2 + y3)) p - x2 (y2 +  y3) q - x3 y3 r;
	DSolve[{equ1 == 0, equ2 == 0, equ3 == 0, equ4 == 0}, f, {x1, x2, x3, y1, y2, y3}]
	\end{lstlisting}
As you see $\textit{DSolve}$ returns un-evaluated i.e. it means that it is not able to solve our system of first order partial differential equations.
\begin{lstlisting}[language=Mathematica,caption={Example code}]
DSolve[2 equ1 + equ2 == 0, f[x1, x2, x3, y1, y2, y3], {x1, x2, x3, y1, y2, y3}][[1, 1]]
(* f[x1, x2, x3, y1, y2, y3] -> C[1][x2/x1, x3/x1, y1, y2, y3] *)
DSolve[equ1 + 2 equ2 == 0, f[x1, x2, x3, y1, y2, y3], {x1, x2, x3, y1, y2, y3}][[1, 1]]
(* f[x1, x2, x3, y1, y2, y3] -> C[1][x1, x2, x3, y2/y1, y3/y1] *)
\end{lstlisting}
Consequently, the dimensionality of this problem can be reduced from six to four.
\begin{lstlisting}[language=Mathematica,caption={Example code}]
f[x1_, x2_, x3_, y1_, y2_, y3_] := g[x2/x1, x3/x1, y2/y1, y3/y1]
equ5 = FullSimplify[(equ3/x1) /. {x2 -> v2 x1, x3 -> v3 x1, y2 -> w2 y1, y3 -> w3 y1}]
(*       (v2 + v3)*w3*Derivative[0, 0, 0, 1][g][v2, v3, w2, w3] + 
v2*w2*Derivative[0, 0, 1, 0][g][v2, v3, w2, w3] - 
v3*(1 + 2*v2 + v3)*Derivative[0, 1, 0, 0][g][v2, v3, w2, w3] - 
v2*(1 + v2)*Derivative[1, 0, 0, 0][g][v2, v3, w2, w3] *)
equ6 = FullSimplify[(equ4/y1) /. {x2 -> v2 x1, x3 -> v3 x1, y2 -> w2 y1, y3 -> w3 y1}] 
(* -(w3*(1 + 2*w2 + w3)*Derivative[0, 0, 0, 1][g][v2, v3, w2, w3]) - 
w2*(1 + w2)*Derivative[0, 0, 1, 0][g][v2, v3, w2, w3] + 
v3*(1 + w2)*Derivative[0, 1, 0, 0][g][v2, v3, w2, w3] + 
v2*Derivative[1, 0, 0, 0][g][v2, v3, w2, w3] *)
\end{lstlisting}
Although $\textit{DSolve}$ cannot solve these equations as a pair either. But  it can solve each separately.
\begin{lstlisting}[language=Mathematica,caption={Example code}]
DSolve[equ5 == 0, g[v2, v3, w2, w3], {v2, v3, w2, w3}][[1, 1]]/.C[1] -> c5 
(* g[v2, v3, w2, w3] -> c5[(v2 (1 + v2 + v3))/v3, (1 + v2) w2, (v3 w3)/v2] *)
(DSolve[equ6 == 0, g[v2, v3, w2, w3], {v2, v3, w2, w3}][[1, 1]] /. 
C[1] -> c6) // FullSimplify 
(* g[v2, v3, w2, w3] -> c6[-((v3 (1 + w2))/v2), ((1 + w2) (1 + w2 + w3))/(v2 w3), 
-Log[(1 + w2)/(v2 w2)]] *)
\end{lstlisting}

The first results indicates that $g$ is a function of
\begin{lstlisting}[language=Mathematica,caption={Example code}]
var5 = List @@ %%[[2]]
(* {(v2 (1 + v2 + v3))/v3, (1 + v2) w2, (v3 w3)/v2} *)
\end{lstlisting}
and also
\begin{lstlisting}[language=Mathematica,caption={Example code}]
var6 = List @@ %%[[2]]
(* {-((v3 (1 + w2))/v2), ((1 + w2) (1 + w2 + w3))/(v2 w3), -Log[(1 + w2)/(v2 w2)]} *)
\end{lstlisting}
The second list of functions can be simplified by
\begin{lstlisting}[language=Mathematica,caption={Example code}]
var6[[3]] = Exp[var6[[3]]]; 
var6[[1]] = -var6[[1]] var6[[3]]; 
var6[[2]] = var6[[2]] var6[[3]]; 
var6 
(* {v3 w2, (w2 (1 + w2 + w3))/w3, (v2 w2)/(1 + w2)} *)
\end{lstlisting}
Now the next step is to combine the previous two expressions for $g$ for to obtain a single expression, presumably as a function of two variables. 

The system of PDEs above can be solved using the procedure described in Chapter V, Sec IV of Goursat's Differential Equations \cite{3}.

The first step is to find the complete, non-commutative group of differential operators that includes $equ5$ and $equ6$.
\begin{lstlisting}[language=Mathematica,caption={Example code}]
comm[equa_, equb_] :=
Collect[(equa /. {Derivative[1, 0, 0, 0][g][v2, v3, w2, w3] -> D[equb, v2], 
Derivative[0, 1, 0, 0][g][v2, v3, w2, w3] -> D[equb, v3], 
Derivative[0, 0, 1, 0][g][v2, v3, w2, w3] -> D[equb, w2], 
Derivative[0, 0, 0, 1][g][v2, v3, w2, w3] -> D[equb, w3]}) - 
(equb /. {Derivative[1, 0, 0, 0][g][v2, v3, w2, w3] -> D[equa, v2], 
Derivative[0, 1, 0, 0][g][v2, v3, w2, w3] -> D[equa, v3], 
Derivative[0, 0, 1, 0][g][v2, v3, w2, w3] -> D[equa, w2], 
Derivative[0, 0, 0, 1][g][v2, v3, w2, w3] -> D[equa, w3]}), 
{Derivative[1, 0, 0, 0][g][v2, v3, w2, w3], Derivative[0, 1, 0, 0][g][v2, v3, w2, w3], 
Derivative[0, 0, 1, 0][g][v2, v3, w2, w3], Derivative[0, 0, 0, 1][g][v2, v3, w2, w3]}, 
Simplify]
equ7 = comm[equ5, equ6]
(* -(w3*(v3*(1 + w2 + w3) + v2*(1 + 2*w2 + w3))*Derivative[0, 0, 0, 1][g][v2, v3, w2, w3])
- v2*w2*(1 + w2)*Derivative[0, 0, 1, 0][g][v2, v3, w2, w3] + v3*(v3*(1 + w2) 
+ v2*(2 + w2))*Derivative[0, 1, 0, 0][g][v2, v3, w2, w3] + 
v2^2*Derivative[1, 0, 0, 0][g][v2, v3, w2, w3] *)
\end{lstlisting}
which by inspection is independent of $equ5$ and $equ6$. On the other hand, $comm[equ5, equ7]$ and $comm[equ6, equ7]$ do not yield independent equations, again by inspection. Thus $\{equ5, equ6, equ7\}$ is a complete group of three operators in four independent variables. From this information alone, we know that $g$ is an arbitrary function of precisely one first integral. This first integral can be obtained by systematically eliminating variables and equations, one pair at a time, until a single equation of two variable remains. We start by solving any one of the equations.
\begin{lstlisting}[language=Mathematica,caption={Example code}]
DSolve[equ5 == 0, g[v2, v3, w2, w3], {v2, v3, w2, w3}][[1, 1]]
(* g[v2, v3, w2, w3] -> C[1][(v2 (1 + v2 + v3))/v3, (1 + v2) w2, (v3 w3)/v2] *)
\end{lstlisting}
and use the solution as the basis for a change of variables:
\vspace*{1cm}

\begin{lstlisting}[language=Mathematica,caption={Example code}]
g[v2_, v3_, w2_, w3_] := h[w2, (v2 (1 + v2 + v3))/v3, (1 + v2) w2, (v3 w3)/v2]
solw2 = equ5/(v2 w2) // Simplify 
(* Derivative[1, 0, 0, 0][h][w2, (v2*(1 + v2 + v3))/v3, (1 + v2)*w2, (v3*w3)/v2] *)
\end{lstlisting}
indicating that $h$ is independent of $w2$. This leaves us with two equations in three variables
\begin{lstlisting}[language=Mathematica,caption={Example code}]
newvar = Solve[Thread[{b1, b2, b3, b4} == List @@ solw2], {v2, v3, w2, w3}] // Flatten; 
(((b3 equ7/(b1 - b3)) // FullSimplify) /. solw2 -> 0 /. newvar) // FullSimplify;\\
Collect[((equ6 // FullSimplify) /. solw2 -> 0 /. newvar) //  
FullSimplify, b1, FullSimplify] + % b1/b3   
equ10 = %% /. Derivative[0, n1_, n2_, n3_][h][b1, b2, b3, b4] ->  
Derivative[n1, n2, n3][h][b2, b3, b4]  
(* b4*(b3 + b4 + b2*b4)*Derivative[0, 0, 1][h][b2, b3, b4] + (1 + b3)* 
(b3*Derivative[0, 1, 0][h][b2, b3, b4] + (1 + b2)*Derivative[1, 0, 0][h][b2, b3, b4]) *) 
equ11 = %% /. Derivative[0, n1_, n2_, n3_][h][b1, b2, b3, b4] ->  
Derivative[n1, n2, n3][h][b2, b3, b4] 
(* (-1 + b4)*b4*Derivative[0, 0, 1][h][b2, b3, b4] + (1 + b2 + b3)
*Derivative[1, 0, 0][h][b2, b3, b4] *)
\end{lstlisting}
Proceeding as before, we next solve one of $equ10$ and $equ11$. (We choose the simpler one.) 
\begin{lstlisting}[language=Mathematica,caption={Example code}]
DSolve[equ11 == 0, h[b2, b3, b4], {b2, b3, b4}][[1, 1, 2]] 
(* h[b2, b3, b4] -> C[1][b3][Log[(1 - b4)/((1 + b2 + b3) b4)] *)
\end{lstlisting}
and use it as the basis for a further change of variables.
\begin{lstlisting}[language=Mathematica,caption={Example code}]
h[b2_, b3_, b4_] := k[b2, b3, (1 - b4)/((1 + b2 + b3) b4)]
solb2 = (equ11/(1 + b2 + b3)) // Simplify 
(* Derivative[1, 0, 0][k][b2, b3, (1 - b4)/(b4 + b2*b4 + b3*b4)] *)
\end{lstlisting}
indicating that $k$ is independent of $b2$. This leaves us with one equation in two variables.
\vspace*{2cm}

\begin{lstlisting}[language=Mathematica,caption={Example code}]
newvar1 = Solve[Thread[{c2, c3, c4} == List @@ solb2], {b2, b3, b4}] // Flatten; 
((equ10 // FullSimplify) /. solb2 -> 0 /. newvar1) // FullSimplify 
equ12 = % /. Derivative[0, n1_, n2_][k][c2, c3, c4] -> Derivative[n1, n2][k][c3, c4] 
(* -((1 + c4 + 2*c3*c4)*Derivative[0, 1][k][c3, c4]) + c3*(1 + c3)
*Derivative[1, 0][k][c3, c4] *)
\end{lstlisting}
Finally, $DSolve$ yields
\begin{lstlisting}[language=Mathematica,caption={Example code}]
DSolve[equ12 == 0, k[c3, c4], {c3, c4}][[1, 1, 2]] 
(* k[c3, c4] -> C[1][c3 (1 + c4 + c3 c4)] *)
\end{lstlisting}
Transforming back to the original independent variables gives
\begin{lstlisting}[language=Mathematica,caption={Example code}]
(((% /. Thread[{c2, c3, c4} -> List @@ solb2]) // Simplify) /. 
Thread[{b1, b2, b3, b4} -> List @@ solw2]) // Simplify
(* C[1][(v2 w2 (1 + (1 + v2) w2 + (1 + v2 + v3) w3))/((v2 + v3 + v3 w2) w3)] *)
\end{lstlisting}

\begin{lstlisting}[language=Mathematica,caption={Example code}]
g[v2_, v3_, w2_, w3_] := %
{equ5, equ6, equ7} // Simplify 
(* {0, 0, 0} *)
\end{lstlisting}
Finally, designating the solution for $g$ as $ansg$,
\begin{lstlisting}[language=Mathematica,caption={Example code}]
(ansg /. {v2 -> x2/x1, v3 -> x3/x1, w2 -> y2/y1, w3 -> y3/y1}) // Simplify
(* C[1][(x2 y2 (x3 y3 + x2 (y2 + y3) + x1 (y1 + y2 + y3)))/(x1 (x2 y1 + x3 (y1 + y2)) y3)] *)
f[x1_, x2_, x3_, y1_, y2_, y3_] := % 
{equ1, equ2, equ3, equ4} 
(* {0, 0, 0, 0} *)
\end{lstlisting}

\subsection{Lattice $W_4$ algebra}.
\vspace*{0.5cm}
\begin{lstlisting}[language=Mathematica,caption={Example code}]
p = D[f[x1, x2, x3, y1, y2, y3, z1, z2, z3], x1];

q = D[f[x1, x2, x3, y1, y2, y3, z1, z2, z3], x2];

r = D[f[x1, x2, x3, y1, y2, y3, z1, z2, z3], x3];

o = D[f[x1, x2, x3, y1, y2, y3, z1, z2, z3], y1];

x = D[f[x1, x2, x3, y1, y2, y3, z1, z2, z3], y2];

a = D[f[x1, x2, x3, y1, y2, y3, z1, z2, z3], y3];

b = D[f[x1, x2, x3, y1, y2, y3, z1, z2, z3], z1];

c = D[f[x1, x2, x3, y1, y2, y3, z1, z2, z3], z2];

d = D[f[x1, x2, x3, y1, y2, y3, z1, z2, z3], z3];

equ1 = 2 x1 p + 2 x2 q + 2 x3 r - y1 o - y2 x - y3 a;

equ2 = -x1 p - x2 q - x3 r - z1 b - z2 c - z3 d + 2 y1 o + 2 y2 x + 2 y3 a;

equ3 = 2 z1 b + 2 z2 c + 2 z3 d - y1 o - y2 x - y3 a;

equ4 = (x1 (x1 + 2 x2 + 2 x3)) p + (x2 (x2 + 
2 x3)) q + (x3^2) r - (y1 (x2 + x3)) o - y2 x3 x;

equ5 = (y1 (y1 + 2 y2 + 2 y3)) o + (y2 (y2 + 2 y3)) x + (y3^2) a - (x1 (y1 + y2 + y3)) p 
- x2 (y2 + y3) q - x3 y3 r  - z1 (y2 + y3) b - z2 y3 c;

equ6 = (z1 (z1 + 2 z2 + 2 z3)) b + (z2 (z2 + 2 z3)) c + (z3^2) d - (y1 (z1 + z2 + z3)) o 
- y2 ( z2 + z3) x - y3 z3 a;
\end{lstlisting}
\begin{lstlisting}[language=Mathematica,caption={Example code}]
DSolve[{HX == 0, HY == 0, , HY == 0, EX == 0, EY == 0, EZ == 0}, 
f[x1, x2, x3, y1, y2, y3, z1, z2, z3], {x1, x2, x3, y1, y2, y3, z1, 
z2, z3}]
\end{lstlisting}
Again $DSolve$ returns un-evaluated, meaning that it can not solve the system of equations.
\vspace*{1cm}
\begin{lstlisting}[language=Mathematica,caption={Example code}]
DSolve[  HX +  3 HZ + 2  HY == 0, 
f[x1, x2, x3, y1, y2, y3, z1, z2, z3], {x1, x2, x3, y1, y2, y3, z1, 
z2, z3}][[1, 1]]
(*f[x1, x2, x3, y1, y2, y3, z1, z2, z3] -> 
C[1][x1, x2, x3, y1, y2, y3, z2/z1, z3/z1]*)
\end{lstlisting}
\begin{lstlisting}[language=Mathematica,caption={Example code}]
DSolve[ HX + HZ + 2 HY == 0, 
f[x1, x2, x3, y1, y2, y3, z1, z2, z3], {x1, x2, x3, y1, y2, y3, z1, 
z2, z3}][[1, 1]]
(* f[x1, x2, x3, y1, y2, y3, z1, z2, z3] -> 
C[1][x1, x2, x3, y2/y1, y3/y1, z1, z2, z3] *)
\end{lstlisting}
\begin{lstlisting}[language=Mathematica,caption={Example code}]
DSolve[ 3  HX +  HZ + 2  HY == 0, 
f[x1, x2, x3, y1, y2, y3, z1, z2, z3], {x1, x2, x3, y1, y2, y3, z1, 
z2, z3}][[1, 1]]
(*f[x1, x2, x3, y1, y2, y3, z1, z2, z3] -> 
C[1][x2/x1, x3/x1, y1, y2, y3, z1, z2, z3]*)
\end{lstlisting}
As before, this computation can be simplified by the substitution,
\begin{lstlisting}[language=Mathematica,caption={Example code}]
f[x1, x2, x3, y1, y2, y3, z1, z2, z3] := g[x2/x1, x3/x1, y2/y1, y3/y1, z2/z1, z3/z1]
\end{lstlisting}
in which case the six equations become
\begin{lstlisting}[language=Mathematica,caption={Example code}]
Simplify[{equ1, equ2, equ3}]
(* {0, 0, 0} *)

equ4 = Simplify[Simplify[equ4]/x1 /. {x2 -> x1 v2, x3 -> x1 v3, y2 -> y1 w2, 
y3 -> y1 w3, z2 -> z1 k2, z3 -> z1 k3}]
(* (v2 + v3)*w3*Derivative[0, 0, 0, 1, 0, 0][g][v2, v3, w2, w3, k2, k3] + 
v2*w2*Derivative[0, 0, 1, 0, 0, 0][g][v2, v3, w2, w3, k2, k3] - 
v3*(1 + 2*v2 + v3)*Derivative[0, 1, 0, 0, 0, 0][g][v2, v3, w2, w3, k2, k3] - 
v2*(1 + v2)*Derivative[1, 0, 0, 0, 0, 0][g][v2, v3, w2, w3, k2, k3] *)

equ5 = Simplify[Simplify[equ5]/y1 /. {x2 -> x1 v2, x3 -> x1 v3, y2 -> y1 w2, 
y3 -> y1 w3, z2 -> z1 k2, z3 -> z1 k3}]
(* k3*(w2 + w3)*Derivative[0, 0, 0, 0, 0, 1][g][v2, v3, w2, w3, k2, k3] + 
k2*w2*Derivative[0, 0, 0, 0, 1, 0][g][v2, v3, w2, w3, k2, k3] - 
w3*(1 + 2*w2 + w3)*Derivative[0, 0, 0, 1, 0, 0][g][v2, v3, w2, w3, k2, k3] - 
w2*(1 + w2)*Derivative[0, 0, 1, 0, 0, 0][g][v2, v3, w2, w3, k2, k3] + 
v3*(1 + w2)*Derivative[0, 1, 0, 0, 0, 0][g][v2, v3, w2, w3, k2, k3] + 
v2*Derivative[1, 0, 0, 0, 0, 0][g][v2, v3, w2, w3, k2, k3] *)

equ6 = Simplify[Simplify[ equ6]/z1 /. {x2 -> x1 v2, x3 -> x1 v3, y2 -> y1 w2, 
y3 -> y1 w3, z2 -> z1 k2, z3 -> z1 k3}]
(* -(k3*(1 + 2*k2 + k3)*Derivative[0, 0, 0, 0, 0, 1][g][v2, v3, w2, w3, k2, k3]) - 
k2*(1 + k2)*Derivative[0, 0, 0, 0, 1, 0][g][v2, v3, w2, w3, k2, k3] + 
(1 + k2)*w3*Derivative[0, 0, 0, 1, 0, 0][g][v2, v3, w2, w3, k2, k3] + 
w2*Derivative[0, 0, 1, 0, 0, 0][g][v2, v3, w2, w3, k2, k3] *)
\end{lstlisting}
As before, this system of first-order $PDEs$ can be solved by using the procedure described in Chapter V, Sec IV of Goursat's Differential Equations.

 The first step is to find the complete, non-commutative group of differential operators that includes $equ4$, $equ5$, and $equ6$. To do so, we use the function $comm$, generalized from $W_3$
\begin{lstlisting}[language=Mathematica,caption={Example code}]
drv = {Derivative[1, 0, 0, 0, 0, 0][g][v2, v3, w2, w3, k2, k3], 
Derivative[0, 1, 0, 0, 0, 0][g][v2, v3, w2, w3, k2, k3], 
Derivative[0, 0, 1, 0, 0, 0][g][v2, v3, w2, w3, k2, k3], 
Derivative[0, 0, 0, 1, 0, 0][g][v2, v3, w2, w3, k2, k3], 
Derivative[0, 0, 0, 0, 1, 0][g][v2, v3, w2, w3, k2, k3], 
Derivative[0, 0, 0, 0, 0, 1][g][v2, v3, w2, w3, k2, k3]};
comm[equa_, equb_] := Collect[
(equa /. {Derivative[1, 0, 0, 0, 0, 0][g][v2, v3, w2, w3, k2, k3] -> D[equb, v2], 
Derivative[0, 1, 0, 0, 0, 0][g][v2, v3, w2, w3, k2, k3] -> D[equb, v3], 
Derivative[0, 0, 1, 0, 0, 0][g][v2, v3, w2, w3, k2, k3] -> D[equb, w2], 
Derivative[0, 0, 0, 1, 0, 0][g][v2, v3, w2, w3, k2, k3] -> D[equb, w3], 
Derivative[0, 0, 0, 0, 1, 0][g][v2, v3, w2, w3, k2, k3] -> D[equb, k2],
Derivative[0, 0, 0, 0, 0, 1][g][v2, v3, w2, w3, k2, k3] -> D[equb, k3]}) - 
(equb /. {Derivative[1, 0, 0, 0, 0, 0][g][v2, v3, w2, w3, k2, k3] -> D[equa, v2], 
Derivative[0, 1, 0, 0, 0, 0][g][v2, v3, w2, w3, k2, k3] -> D[equa, v3], 
Derivative[0, 0, 1, 0, 0, 0][g][v2, v3, w2, w3, k2, k3] -> D[equa, w2], 
Derivative[0, 0, 0, 1, 0, 0][g][v2, v3, w2, w3, k2, k3] -> D[equa, w3], 
Derivative[0, 0, 0, 0, 1, 0][g][v2, v3, w2, w3, k2, k3] -> D[equa, k2], 
Derivative[0, 0, 0, 0, 0, 1][g][v2, v3, w2, w3, k2, k3] -> D[equa, k3]}),
drv, Simplify]

equ7 = comm[equ4, equ5]
(* (k3*v2*w2 + k3*(v2 + v3)*w3)*Derivative[0, 0, 0, 0, 0, 1][g][v2, v3, w2, w3, k2, k3] +
k2*v2*w2*Derivative[0, 0, 0, 0, 1, 0][g][v2, v3, w2, w3, k2, k3] - 
w3*(v3*(1 + w2 + w3) + v2*(1 + 2*w2 + w3))*
Derivative[0, 0, 0, 1, 0, 0][g][v2, v3, w2, w3, k2, k3] - 
v2*w2*(1 + w2)*Derivative[0, 0, 1, 0, 0, 0][g][v2, v3, w2, w3, k2, k3] + 
v3*(v3*(1 + w2) + v2*(2 + w2))*Derivative[0, 1, 0, 0, 0, 0][g][v2, v3, w2, w3, k2, k3] + 
v2^2*Derivative[1, 0, 0, 0, 0, 0][g][v2, v3, w2, w3, k2, k3] *)

equ8 = comm[equ5, equ6]
(* -(k3*((1 + 2*k2 + k3)*w2 + (1 + k2 + k3)*w3)*
Derivative[0, 0, 0, 0, 0, 1][g][v2, v3, w2, w3, k2, k3]) - 
k2*(1 + k2)*w2*Derivative[0, 0, 0, 0, 1, 0][g][v2, v3, w2, w3, k2, k3] + 
w3*((2 + k2)*w2 + (1 + k2)*w3)*Derivative[0, 0, 0, 1, 0, 0][g][v2, v3, w2, w3, k2, k3] +
w2^2*Derivative[0, 0, 1, 0, 0, 0][g][v2, v3, w2, w3, k2, k3] - 
v3*w2*Derivative[0, 1, 0, 0, 0, 0][g][v2, v3, w2, w3, k2, k3] *)
\end{lstlisting}
which are independent of the first three operators, increasing the size of the group to five. $comm[equ4, equ6]$ vanishes identically and so does not add an operator. On the other hand, the seven additional commutators involving $equ7$ and $equ8$ yield expressions that are linear combinations of $\{equ4, equ5, equ6, equ7, equ8\}$. Thus, these five operators comprise the entire group.

From this information alone, we know that $g$ is an arbitrary function of precisely one first integral. This first integral can be obtained by systematically eliminating variables and equations, one pair at a time, until a single equation of two variable remains. Start by solving any one of the equations. 
\begin{lstlisting}[language=Mathematica,caption={Example code}]
DSolve[equ4 == 0, 
g[v2, v3, w2, w3, k2, k3], {v2, v3, w2, w3, k2, k3}][[1, 
1]] // FullSimplify
(* g[v2, v3, w2, w3, k2, k3] -> 
C[1][k2, k3][(v2 (1 + v2 + v3))/v3, (1 + v2) w2, (v3 w3)/v2] *)
\end{lstlisting}
\begin{lstlisting}[language=Mathematica,caption={Example code}]
g[v2_, v3_, w2_, w3_, k2_, k3_] := 
h[w2, (v2 (1 + v2 + v3))/v3, (1 + v2) w2, (v3 w3)/v2, k2, k3];
tr1 = {equ4, equ5, equ6, equ7, equ8} // Simplify;
\end{lstlisting}
\begin{lstlisting}[language=Mathematica,caption={Example code}]
solw2 = equ4/(v2 w2) // FullSimplify;
\end{lstlisting}
\begin{lstlisting}[language=Mathematica,caption={Example code}]
newvar = Solve[
Thread[{b1, b2, b3, b4, b5, b6} == List @@ solw2], {v2, v3, w2, 
w3, k2, k3}] // Flatten;
tr1p = Collect[FullSimplify[Rest[tr1] /. solw2 -> 0 /. newvar], b1, 
FullSimplify] /. b1*(z_) -> 0 /. 
Derivative[0, n1_, n2_, n3_, n4_, n5_][h][b1, b2, b3, b4, b5, 
b6] -> Derivative[n1, n2, n3, n4, n5][h][b2, b3, b4, b5, b6];
\end{lstlisting}
\begin{lstlisting}[language=Mathematica,caption={Example code}]
DSolve[First@tr1p == 0, 
h[b2, b3, b4, b5, b6], {b2, b3, b4, b5, b6}] /. Log[z_] -> z;
h[b2_, b3_, b4_, b5_, b6_] := 
j[b4, b3, b5, (1 + b2 + b3) b4 b6, b6 (1 - b4)];
tr2 = tr1p // Simplify;
solb4 = First@tr2/(b4 (b4 - 1)) // Simplify;
newvar = Solve[
Thread[{c1, c2, c3, c4, c5} == List @@ solb4], {b2, b3, b4, b5, 
b6}] // Flatten;

tr2p = Collect[(Cancel[(c1 - 1) Rest@tr2] /. solb4 -> 0 /. newvar) // 
FullSimplify, c1, FullSimplify];
tr2p[[3]] = tr2p[[3]]/c1;
tr2p = tr2p /. c1 z__ -> 0 /. 
Derivative[0, n1_, n2_, n3_, n4_][j][c1, c2, c3, c4, c5]   -> 
Derivative[n1, n2, n3, n4][j][c2, c3, c4, c5];
\end{lstlisting}
\begin{lstlisting}[language=Mathematica,caption={Example code}]
DSolve[Last@tr2p == 0, j[c2, c3, c4, c5], {c2, c3, c4, c5}];
j[c2_, c3_, c4_, c5_] := l[c5, c2, c3, (c2 - c4)/(1 + c3 + c5)];
tr3 = -tr2p // Simplify // RotateRight;
solc5 = First@tr3/(c5 (1 + c3 + c5)) // Simplify;
newvar = Solve[
Thread[{d1, d2, d3, d4} == List @@ solc5], {c2, c3, c4, c5}] // 
Flatten;
tr3p = Collect[(Rest@tr3 /. solc5 -> 0 /. newvar) // FullSimplify, d1,
FullSimplify] /. d1 z_ -> 0 /. 
Derivative[0, n1_, n2_, n3_][l][d1, d2, d3, d4] -> 
Derivative[n1, n2, n3][l][d2, d3, d4];
\end{lstlisting}
\begin{lstlisting}[language=Mathematica,caption={Example code}]
DSolve[Last@tr3p == 0, l[d2, d3, d4], {d2, d3, d4}] // Simplify;
l[d2_, d3_, d4_] := m[d3, (1 + d2) d3, (d2 (1 + d2 - d4))/d4];
tr4 = tr3p // Simplify // RotateRight;
sold3 = First@tr4/(d2 d3);
newvar = Solve[Thread[{e1, e2, e3} == List @@ sold3], {d2, d3, d4}] //
Flatten;
tr4p = Collect[(-(e2/e1) Rest@tr4 /. sold3 -> 0 /. newvar) // 
FullSimplify, e1, FullSimplify] /. e1 z_ -> 0 /. 
Derivative[0, n1_, n2_][m][e1, e2, e3] -> 
Derivative[n1, n2][m][e2, e3];
\end{lstlisting}
\begin{lstlisting}[language=Mathematica,caption={Example code}]
(DSolve[Last@tr4p == 0, m[e2, e3], {e2, e3}] // Simplify) /. 
Log[z_] -> z;
((((((((%[[1, 1, 2]] /. Thread[{e1, e2, e3} -> List @@ sold3]) // 
Simplify) /. 
Thread[{d1, d2, d3, d4} -> List @@ solc5]) // 
Simplify) /. 
Thread[{c1, c2, c3, c4, c5} -> List @@ solb4]) // 
Simplify) /. 
Thread[{b1, b2, b3, b4, b5, b6} -> List @@ solw2] // 
Simplify) /. {v2 -> x2/x1, v3 -> x3/x1, w2 -> y2/y1, w3 -> y3/y1,
k2 -> z2/z1, k3 -> z3/z1})  // Simplify
\end{lstlisting}
Final solution
\begin{lstlisting}[language=Mathematica,caption={Example code}]
C[1][-(((x2 y2 z2 + x1 (y2 z2 + y1 (z1 + z2))) (x3 y3 z3 + 
x2 (y3 z3 + y2 (z2 + z3))))/(
x2 y2 z2 (x3 y3 z3 + x2 (y3 z3 + y2 (z2 + z3)) + 
x1 (y3 z3 + y2 (z2 + z3) + y1 (z1 + z2 + z3)))))]
\end{lstlisting}

\subsection{Expressing a fractional multivariate polynomials to its low-order polynomial decomposition}

Suppose we have given the following question.

\textbf{Question:}

Let $f2$ be fractional multivariate polynomial as follows
\begin{lstlisting}[language=Mathematica,caption={Example code}]
f2 = -((2 x1 x2 x3 x4 y1 y2^2 y3 (x2 y1 + (x3 + x4) (y1 + y2) + x4 y3) (x3 y3 + x2 (y2 + y3) + x1 (y1 + y2 + y3)))/((x2 y2 + x1 (y1 + y2))^2 (x2 y1 + x3 (y1 + y2)) (x3 y3 + x2 (y2 + y3)) (x3 y2 + x4 (y2 + y3))^2));

\end{lstlisting}
and also let $k1$ and $k2$ be given as follows

\begin{lstlisting}[language=Mathematica,caption={Example code}]
k1 = (x2 y2 (x3 y3 + x2 (y2 + y3) + x1 (y1 + y2 + y3)))/((x2 y2 + x1 (y1 + y2)) (x3 y3 + x2 (y2 + y3)));
k2 = (x3 y2 (x2 y1 + (x3 + x4) (y1 + y2) + x4 y3))/((x2 y1 + x3 (y1 + y2)) (x3 y2 + x4 (y2 + y3)));
\end{lstlisting}

then express $f2$ as a low - order polynomial in $k1$ and $k2$.\\

This can be done as follows. \\
First, generate a generic low order polynomial.
\begin{lstlisting}[language=Mathematica,caption={Example code}]
Map[t1^First@# t2^Last@# &, Tuples[Range[0, 3], 2]].Table[Unique["c"], {16}]
(* c3 + c7 t1 + c11 t1^2 + c15 t1^3 + c4 t2 + c8 t1 t2 + c12 t1^2 t2 + 
c16 t1^3 t2 + c5 t2^2 + c9 t1 t2^2 + c13 t1^2 t2^2 + c17 t1^3 t2^2 + 
c6 t2^3 + c10 t1 t2^3 + c14 t1^2 t2^3 + c18 t1^3 t2^3 *)
\end{lstlisting}
and then use SolveAlways. After about twenty seconds we will gwt result
\begin{lstlisting}[language=Mathematica,caption={Example code}]
Flatten@SolveAlways[f2 == (% /. {t1 -> k1, t2 -> k2}), {x1, x2, x3, x4, y1, y2, y3}]
(* {c3 -> 0, c4 -> 0, c5 -> 0, c6 -> 0, c11 -> 0, c15 -> 0, c7 -> 0, c12 -> 2, c16 -> 0, c8 -> -2, c10 -> 0, c13 -> -2, c14 -> 0, c17 -> 0, c18 -> 0, c9 -> 2} *)
\end{lstlisting}
And we have the final solution
\begin{lstlisting}[language=Mathematica,caption={Example code}]
Factor[%% /. %]
(* -2 (-1 + t1) t1 (-1 + t2) t2 *)
\end{lstlisting}
which is the desired result.

And for completeness we have
\begin{lstlisting}[language=Mathematica,caption={Example code}]
Simplify[f2 == % /. {t1 -> k1, t2 -> k2}]
(* True *)
\end{lstlisting}
\textbf{Also here we have much faster alternative:}

Because SolveAlways determines the coefficients $c$ for any $\{x1, x2, x3, x4, y1, \\
y2, y3\}$, Solve must be able to obtain the same values for the coefficients $c$ for specific values of $\{x1, x2, x3, x4, y1, y2, y3\}$, and much faster.
As before we do have $f2$ and $k1$ and $k2$.
\begin{lstlisting}[language=Mathematica,caption={Example code}]
f2 = -((2 x1 x2 x3 x4 y1 y2^2 y3 (x2 y1 + (x3 + x4) (y1 + y2) + x4 y3) (x3 y3 + x2 (y2 + y3) + x1 (y1 + y2 + y3)))/((x2 y2 + x1 (y1 + y2))^2 (x2 y1 + x3 (y1 + y2)) (x3 y3 + x2 (y2 + y3)) (x3 y2 + x4 (y2 + y3))^2));

\end{lstlisting}
\begin{lstlisting}[language=Mathematica,caption={Example code}]
k1 = (x2 y2 (x3 y3 + x2 (y2 + y3) + x1 (y1 + y2 + y3)))/((x2 y2 + x1 (y1 + y2)) (x3 y3 + x2 (y2 + y3)));
k2 = (x3 y2 (x2 y1 + (x3 + x4) (y1 + y2) + x4 y3))/((x2 y1 + x3 (y1 + y2)) (x3 y2 + x4 (y2 + y3)));
\end{lstlisting}
\begin{lstlisting}[language=Mathematica,caption={Example code}]
tp = Tuples[Range[0, 3], 2]; tp // Length
(* 16 *)
\end{lstlisting}
\begin{lstlisting}[language=Mathematica,caption={Example code}]
gp = Map[t1^#[[1]] t2^#[[2]] &, tp].Table[Unique["c"], {tp // Length}]
(* c3 + c7 t1 + c11 t1^2 + c15 t1^3 + c4 t2 + c8 t1 t2 + c12 t1^2 t2 + c16 t1^3 t2 + c5 t2^2 + c9 t1 t2^2 + c13 t1^2 t2^2 + c17 t1^3 t2^2 + c6 t2^3 + c10 t1 t2^3 + c14 t1^2 t2^3 + c18 t1^3 t2^3 *)
\end{lstlisting}
\begin{lstlisting}[language=Mathematica,caption={Example code}]
Flatten@Solve[Table[(f2 == (gp /. {t1 -> k1, t2 -> k2})) /. Thread[{x1, x2, x3, x4, y1, y2, y3} -> RandomInteger[{1, 7}, 7]], {n, tp // Length}], List @@ (First@# & /@ (gp /. gp[[1]] -> gp[[1]] z))]
(* {c10 -> 0, c11 -> 0, c12 -> 2, c13 -> -2, c14 -> 0, c15 -> 0, c16 -> 0, c17 -> 0, c18 -> 0, c3 -> 0, c4 -> 0, c5 -> 0, c6 -> 0, c7 -> 0, c8 -> -2, c9 -> 2} *)
\end{lstlisting}
\begin{lstlisting}[language=Mathematica,caption={Example code}]
Factor[%% /. %]
(* -2 (-1 + t1) t1 (-1 + t2) t2 *)
\end{lstlisting}
\begin{lstlisting}[language=Mathematica,caption={Example code}]
Simplify[f2 == % /. {t1 -> k1, t2 -> k2}]
(* True *)
\end{lstlisting}

\textbf{Question:}

Let $f6$ be fractional multivariate polynomial as follows
\begin{lstlisting}[language=Mathematica,caption={Example code}]
f6 = (2 x1 x2 x5 x6 y2 (x2 y1 + x3 (y1 + y2)) y3^2 y4 (x5 y5 + x4 (y4 + y5)))/((x2 y2 + x1 (y1 + y2)) (x3 y3 + x2 (y2 + y3))^2 (x4 y3 + x5 (y3 + y4))^2 (x5 y4 + x6 (y4 + y5)));

\end{lstlisting}
and also let $k1$, $k2$, $k3$, $k4$, $k5$ and $k6$ be given as follows

\begin{lstlisting}[language=Mathematica,caption={Example code}]
k1 = ((x2 y2 + x1 (y1 + y2)) (x3 y3 + x2 (y2 + y3)))/(x2 y2 (x3 y3 + x2 (y2 + y3) + x1 (y1 + y2 + y3)));
k2 = ((x2 y1 + x3 (y1 + y2)) (x3 y2 +  x4 (y2 + y3)))/(x3 y2 (x2 y1 + (x3 + x4) (y1 + y2) + x4 y3));
k3 = ((x3 y3 + x2 (y2 + y3)) (x4 y4 + x3 (y3 + y4)))/(x3 y3 (x4 y4 +  x3 (y3 + y4) + x2 (y2 + y3 + y4)));
k4 = ((x3 y2 + x4 (y2 + y3)) (x4 y3 + x5 (y3 + y4)))/(x4 y3 (x3 y2 + (x4 + x5) (y2 + y3) + x5 y4));
k5 = ((x4 y4 + x3 (y3 + y4)) (x5 y5 + x4 (y4 + y5)))/(x4 y4 (x5 y5 + x4 (y4 + y5) + x3 (y3 + y4 + y5)));
k6 = ((x4 y3 + x5 (y3 + y4)) (x5 y4 + x6 (y4 + y5)))/(x5 y4 (x4 y3 + (x5 + x6) (y3 + y4) + x6 y5));
\end{lstlisting}

then express $f6$ as a low - order polynomial in $k1$, $k2$, $k3$, $k4$, $k5$ and $k6$.

\begin{lstlisting}[language=Mathematica,caption={Example code}]
tp = Tuples[Range[-1, 1], 6]; tp // Length
(* 729 *)
\end{lstlisting}

\begin{lstlisting}[language=Mathematica,caption={Example code}]
gp = Map[t1^#[[1]] t2^#[[2]] t3^#[[3]] t4^#[[4]] t5^#[[5]] t6^#[[6]] &, tp].Table[Unique["c"], {tp // Length}];
\end{lstlisting}

\begin{lstlisting}[language=Mathematica,caption={Example code}]
sol = Flatten@ Solve[Table[(f6 == (gp /. {t1 -> k1, t2 -> k2, t3 -> k3, t4 -> k4, t5 -> k5, t6 -> k6})) /. Thread[{x1, x2, x3, x4, x5, x6, y1, y2, y3, y4, y5} -> RandomInteger[{1, 11}, 11]], {n, tp // Length}], List @@ (First@# & /@ (gp /. gp[[1]] -> gp[[1]] z))]; sol /. Rule[_, 0] -> Nothing
(* {Nothing, ..., Nothing, c114 -> -2, c115 -> 2, Nothing,..., Nothing, c123 -> 2, c124 -> -2, Nothing, ..., Nothing, c330 -> -2, c331 -> 2, Nothing, ... , Nothing, c339 -> 2, Nothing, c340 -> -2, Nothing, ..., Nothing, c357 -> 2, c358 -> -2, Nothing, ..., Nothing, c366 -> -2, c367 -> 2, Nothing, ..., Nothing, c87 -> 2, c88 -> -2, Nothing, ..., Nothing, c96 -> -2, c97 -> 2, Nothing, Nothing} *)
\end{lstlisting}

\begin{lstlisting}[language=Mathematica,caption={Example code}]
Factor[gp /. sol]
(* (2 (-1 + t1) (-1 + t3) (-1 + t4) (-1 + t6))/(t1 t3 t4 t6) *)
\end{lstlisting}
Which is the desired result.

And as before for completeness we have
\begin{lstlisting}[language=Mathematica,caption={Example code}]
Simplify[f6 == % /. {t1 -> k1, t2 -> k2, t3 -> k3, t4 -> k4, t5 -> k5, t6 -> k6}]
(* True *)
\end{lstlisting}

\textbf{By using Groebner Basis:}

Also there is another way for to reach to the solution by using Groebner Basis. But this approach is very slow!

\begin{lstlisting}[language=Mathematica,caption={Example code}]
poly = (2 x1 x2 x5 x6 y2 (x2 y1 + x3 (y1 + y2)) y3^2 y4 (x5 y5 + x4 (y4 + y5)))/((x2 y2 + x1 (y1 + y2)) (x3 y3 + x2 (y2 + y3))^2 (x4 y3 +  x5 (y3 + y4))^2 (x5 y4 + x6 (y4 + y5)));
\end{lstlisting}

\begin{lstlisting}[language=Mathematica,caption={Example code}]
eqns = {K1 == ((x2 y2 + x1 (y1 + y2)) (x3 y3 + x2 (y2 + y3)))/(x2 y2 (x3 y3 + x2 (y2 + y3) + x1 (y1 + y2 + y3))), 
K2 == ((x2 y1 + x3 (y1 + y2)) (x3 y2 +  x4 (y2 + y3)))/(x3 y2 (x2 y1 + (x3 + x4) (y1 + y2) + x4 y3)),
K3 == ((x3 y3 + x2 (y2 + y3)) (x4 y4 +  x3 (y3 + y4)))/(x3 y3 (x4 y4 + x3 (y3 + y4) + x2 (y2 + y3 + y4))), 
K4 == ((x3 y2 + x4 (y2 + y3)) (x4 y3 +  x5 (y3 + y4)))/(x4 y3 (x3 y2 + (x4 + x5) (y2 + y3) + x5 y4)),
K5 == ((x4 y4 + x3 (y3 + y4)) (x5 y5 +  x4 (y4 + y5)))/(x4 y4 (x5 y5 + x4 (y4 + y5) +  x3 (y3 + y4 + y5))), 
K6 == ((x4 y3 + x5 (y3 + y4)) (x5 y4 +  x6 (y4 + y5)))/(x5 y4 (x4 y3 + (x5 + x6) (y3 + y4) + x6 y5))};
\end{lstlisting}

Now let us compute Groebner Basis

\begin{lstlisting}[language=Mathematica,caption={Example code}]
gb = GroebnerBasis[eqns, {x1, y1, x2, y2, x3, y3, x4, y4, x5, y5, x6}];
\end{lstlisting}

The remainder $r$ gives a representation of poly in terms of $K1$ , $K2$, $K3$, $K4$, $K5$ and $ K6$.

\begin{lstlisting}[language=Mathematica,caption={Example code}]
{qs, r} = PolynomialReduce[poly, gb, {x1, y1, x2, y2, x3, y3, x4, y4, x5, y5, x6}];
\end{lstlisting}

Where $r$ is our solution in $K1$ , $K2$, $K3$, $K4$, $K5$ and $ K6$. And the following code validates correctness:

\begin{lstlisting}[language=Mathematica,caption={Example code}]
poly == r /. ToRules[And @@ eqns] // Expand
\end{lstlisting}
And please note that, this may take a while. (May be more than a while! It depends on how powerful is your computer. )

\subsection{Checking symmetries in our shift operators}. 
First, before starting, we need to know which variables are employed in our functions. For to do this we employ the following code:

Set 

\begin{lstlisting}[language=Mathematica,caption={Example code}]
f9 = (2 x1 x2 x5 y2 y5 y6 z2 (x2 y1 y2 z1 + x2 y1 y3 z1 + x3 y1 y3 z1 + x2 y1 y3 z2 + x3 y1 y3 z2 + x3 y2 y3 z2) z3^2 z4 (x4 x5 y4 z4 + x4 x6 y4 z4 + x4 x6 y5 z4 + x4 x6 y4 z5 + x4 x6 y5 z5 + x5 x6 y5 z5))/((x1 y1 z1 + x1 y1 z2 + x1 y2 z2 + x2 y2 z2) (x2 y2 z2 + x2 y2 z3 + x2 y3 z3 +  x3 y3 z3)^2 (x4 y4 z3 + x4 y5 z3 + x5 y5 z3 + x5 y5 z4)^2 (x5 y5 z4 + x5 y6 z4 + x6 y6 z4 + x6 y6 z5));
\end{lstlisting}
and
\begin{lstlisting}[language=Mathematica,caption={Example code}]
k1 = ((x1 y1 z1 + x2 y2 z2 + x1 (y1 + y2) z2) (x2 y2 z2 +  x3 y3 z3 + x2 (y2 + y3) z3))/(x2 y2 z2 (x2 y2 z2 + x3 y3 z3 + x2 (y2 + y3) z3 + x1 (y2 z2 + (y2 + y3) z3 + y1 (z1 + z2 + z3)))); 
k2 = (((x2 + x3) y1 z1 + x3 (y1 + y2) z2) ((x3 + x4) y2 z2 + x4 (y2 + y3) z3))/(x3 y2 z2 (x2 y1 z1 + (x3 + x4) (y2 z2 + y1 (z1 + z2)) + x4 (y1 + y2 + y3) z3)); 
k3 = ((x2 (y2 + y3) z1 + x3 y3 (z1 + z2)) (x3 (y3 + y4) z2 + x4 y4 (z2 + z3)))/(x3 y3 z2 (x2 (y2 + y3 + y4) z1 + x3 (y3 + y4) (z1 + z2) + x4 y4 (z1 + z2 + z3))); 
k4 = ((x2 y2 z2 + x3 y3 z3 + x2 (y2 + y3) z3) (x3 y3 z3 +  x4 y4 z4 + x3 (y3 + y4) z4))/(x3 y3 z3 (x3 y3 z3 + x4 y4 z4 +  x3 (y3 + y4) z4 + x2 (y3 z3 + (y3 + y4) z4 + y2 (z2 + z3 + z4)))); 
k5 = (((x3 + x4) y2 z2 + x4 (y2 + y3) z3) ((x4 + x5) y3 z3 + x5 (y3 + y4) z4))/(x4 y3 z3 (x3 y2 z2 + (x4 + x5) (y3 z3 + y2 (z2 + z3)) + x5 (y2 + y3 + y4) z4)); 
k6 = ((x3 (y3 + y4) z2 + x4 y4 (z2 + z3)) (x4 (y4 + y5) z3 +  x5 y5 (z3 + z4)))/(x4 y4 z3 (x3 (y3 + y4 + y5) z2 + x4 (y4 + y5) (z2 + z3) + x5 y5 (z2 + z3 + z4))); 
k7 = ((x3 y3 z3 + x4 y4 z4 + x3 (y3 + y4) z4) (x4 y4 z4 + x5 y5 z5 + x4 (y4 + y5) z5))/(x4 y4 z4 (x4 y4 z4 + x5 y5 z5 + x4 (y4 + y5) z5 + x3 (y4 z4 + (y4 + y5) z5 + y3 (z3 + z4 + z5)))); 
k8 = (((x4 + x5) y3 z3 + x5 (y3 + y4) z4) ((x5 + x6) y4 z4 + x6 (y4 + y5) z5))/(x5 y4 z4 (x4 y3 z3 + (x5 + x6) (y4 z4 + y3 (z3 + z4)) + x6 (y3 + y4 + y5) z5)); 
k9 = ((x4 (y4 + y5) z3 + x5 y5 (z3 + z4)) (x5 (y5 + y6) z4 + x6 y6 (z4 + z5)))/(x5 y5 z4 (x4 (y4 + y5 + y6) z3 + x5 (y5 + y6) (z3 + z4) + x6 y6 (z3 + z4 + z5)));
\end{lstlisting}
Then by using the following code we will obtain the set of our variables which have been employed
\begin{lstlisting}[language=Mathematica,caption={Example code}]
Union[Cases[#, _Symbol, Infinity]] & /@ {f9}
(* {{x1, x2, x3, x4, x5, x6, y1, y2, y3, y4, y5, y6, z1, z2, z3, z4, z5}} *)
\end{lstlisting}
\begin{lstlisting}[language=Mathematica,caption={Example code}]
Union[Cases[#, _Symbol, Infinity]] & /@ {k1, k2, k3, k4, k5, k6, k7, k8, k9}
(* {{x1, x2, x3, y1, y2, y3, z1, z2, z3}, {x2, x3, x4, y1, y2, y3, z1, z2, z3}, {x2, x3, x4, y2, y3, y4, z1, z2, z3}, {x2, x3, x4, y2, y3, y4, z2, z3, z4}, {x3, x4, x5, y2, y3, y4, z2, z3, z4}, {x3, x4, x5, y3, y4, y5, z2, z3, z4}, {x3, x4, x5, y3, y4, y5, z3, z4, z5}, {x4, x5, x6, y3, y4, y5, z3, z4, z5}, {x4, x5, x6, y4, y5, y6, z3, z4, z5}} *)
\end{lstlisting}

Now in what comes below, we  specifically mean that for example in $sl_4$ in a process for finding $F_{9}^{(4)}$, the substitution 
$$\{x1, x2, x3, x4, x5, x6, y1, y2, y3, y4, y5, y6, z1, z2, z3, z4, z5 \} $$ 
instead of 
$$ \{y6, y5, y4, y3, y2, y1, x6, x5, x4, x3, x2, x1, z5, z4, z3, z2, z1\}$$
transforms $\tau_{9}^{(4)}$ to $\tau_{1}^{(4)}$, $\tau_{8}^{(4)}$ to $\tau_{2}^{(4)}$, $\tau_{7}^{(4)}$ to $\tau_{3}^{(4)}$, and $\tau_{6}^{(4)}$ to $\tau_{4}^{(4)}$ while leaving $F_{9}^{(4)}$ unchanged.

Therefore, those four pairs must enter the expression for $F_{9}^{(4)}$ symmetrically. 

As a result, the generic polynomials we have been using above, can be reduced greatly in numbers of terms, a factor of $(\frac{2}{3})^4$.  Corresponding running time then should be reduced by a factor of $(\frac{2}{3})^8$, other things being equal.  It is possible that additional symmetries exist! It needs to be checked!

Here for simplification and for to be able in codding them, we write instead $F_{9}^{(4)} := F9 $ and $\tau_{7}^{(4)} := Ki $;

\begin{lstlisting}[language=Mathematica,caption={Example code}]
arg1 = {a1, a2, a3, a4, a5, a6, b1, b2, b3, b4, b5, b6, d1, d2, d3, d4, d5};
arg2 = {a6, a5, a4, a3, a2, a1, b6, b5, b4, b3, b2, b1, d5, d4, d3, d2, d1};
arg3 = {b6, b5, b4, b3, b2, b1, a6, a5, a4, a3, a2, a1, d5, d4, d3, d2, d1};
\end{lstlisting}

\begin{lstlisting}[language=Mathematica,caption={Example code}]
F9[x1_, x2_, x3_, x4_, x5_, x6_, y1_, y2_, y3_, y4_, y5_, y6_, z1_, z2_, z3_, z4_, z5_ ] := (2 x1 x2 x5 y2 y5 y6 z2 (x2 y1 y2 z1 + x2 y1 y3 z1 + 
x3 y1 y3 z1 + x2 y1 y3 z2 + x3 y1 y3 z2 + x3 y2 y3 z2) z3^2 z4 (x4 x5 y4 z4 + x4 x6 y4 z4 + x4 x6 y5 z4 + x4 x6 y4 z5 + x4 x6 y5 z5 + x5 x6 y5 z5))/((x1 y1 z1 + x1 y1 z2 + x1 y2 z2 + x2 y2 z2) (x2 y2 z2 + x2 y2 z3 + x2 y3 z3 + x3 y3 z3)^2 (x4 y4 z3 + x4 y5 z3 + x5 y5 z3 + x5 y5 z4)^2 (x5 y5 z4 + x5 y6 z4 + x6 y6 z4 + x6 y6 z5))
\end{lstlisting}

\begin{lstlisting}[language=Mathematica,caption={Example code}]
Simplify[F9 @@ arg1 == F9 @@ arg3]
(* True *)
\end{lstlisting}

\begin{lstlisting}[language=Mathematica,caption={Example code}]
K1[x1_, x2_, x3_, x4_, x5_, x6_, y1_, y2_, y3_, y4_, y5_, y6_, z1_, z2_, z3_, z4_, z5_] := ((x1 y1 z1 + x2 y2 z2 + x1 (y1 + y2) z2) (x2 y2 z2 +  x3 y3 z3 + x2 (y2 + y3) z3))/(x2 y2 z2 (x2 y2 z2 + x3 y3 z3 + x2 (y2 + y3) z3 + x1 (y2 z2 + (y2 + y3) z3 + y1 (z1 + z2 + z3)))); 
K2[x1_, x2_, x3_, x4_, x5_, x6_, y1_, y2_, y3_, y4_, y5_, y6_, z1_, z2_, z3_, z4_, z5_] := (((x2 + x3) y1 z1 + x3 (y1 + y2) z2) ((x3 + x4) y2 z2 + x4 (y2 + y3) z3))/(x3 y2 z2 (x2 y1 z1 + (x3 + x4) (y2 z2 + y1 (z1 + z2)) + x4 (y1 + y2 + y3) z3)); 
K3[x1_, x2_, x3_, x4_, x5_, x6_, y1_, y2_, y3_, y4_, y5_, y6_, z1_, z2_, z3_, z4_, z5_] := ((x2 (y2 + y3) z1 + x3 y3 (z1 + z2)) (x3 (y3 + y4) z2 + x4 y4 (z2 + z3)))/(x3 y3 z2 (x2 (y2 + y3 + y4) z1 + x3 (y3 + y4) (z1 + z2) + x4 y4 (z1 + z2 + z3))); 
K4[x1_, x2_, x3_, x4_, x5_, x6_, y1_, y2_, y3_, y4_, y5_, y6_, z1_, z2_, z3_, z4_,  z5_] := ((x2 y2 z2 + x3 y3 z3 + x2 (y2 + y3) z3) (x3 y3 z3 +  x4 y4 z4 + x3 (y3 + y4) z4))/(x3 y3 z3 (x3 y3 z3 + x4 y4 z4 +  x3 (y3 + y4) z4 + x2 (y3 z3 + (y3 + y4) z4 + y2 (z2 + z3 + z4)))); 
K5[x1_, x2_, x3_, x4_, x5_, x6_, y1_, y2_, y3_, y4_, y5_, y6_, z1_, z2_, z3_, z4_, z5_] := (((x3 + x4) y2 z2 + x4 (y2 + y3) z3) ((x4 + x5) y3 z3 + x5 (y3 + y4) z4))/(x4 y3 z3 (x3 y2 z2 + (x4 + x5) (y3 z3 + y2 (z2 + z3)) + x5 (y2 + y3 + y4) z4)); 
K6[x1_, x2_, x3_, x4_, x5_, x6_, y1_, y2_, y3_, y4_, y5_, y6_, z1_, z2_, z3_, z4_, z5_] := ((x3 (y3 + y4) z2 + x4 y4 (z2 + z3)) (x4 (y4 + y5) z3 +  x5 y5 (z3 + z4)))/(x4 y4 z3 (x3 (y3 + y4 + y5) z2 + x4 (y4 + y5) (z2 + z3) + x5 y5 (z2 + z3 + z4))); 
K7[x1_, x2_, x3_, x4_, x5_, x6_, y1_, y2_, y3_, y4_, y5_, y6_, z1_, z2_, z3_, z4_, z5_] := ((x3 y3 z3 + x4 y4 z4 + x3 (y3 + y4) z4) (x4 y4 z4 + x5 y5 z5 + x4 (y4 + y5) z5))/(x4 y4 z4 (x4 y4 z4 + x5 y5 z5 + x4 (y4 + y5) z5 + x3 (y4 z4 + (y4 + y5) z5 + y3 (z3 + z4 + z5)))); 
K8[x1_, x2_, x3_, x4_, x5_, x6_, y1_, y2_, y3_, y4_, y5_, y6_, z1_, z2_, z3_, z4_, z5_] := (((x4 + x5) y3 z3 + x5 (y3 + y4) z4) ((x5 + x6) y4 z4 + x6 (y4 + y5) z5))/(x5 y4 z4 (x4 y3 z3 + (x5 + x6) (y4 z4 + y3 (z3 + z4)) + x6 (y3 + y4 + y5) z5)); 
K9[x1_, x2_, x3_, x4_, x5_, x6_, y1_, y2_, y3_, y4_, y5_, y6_, z1_, z2_, z3_, z4_, z5_] := ((x4 (y4 + y5) z3 + x5 y5 (z3 + z4)) (x5 (y5 + y6) z4 + x6 y6 (z4 + z5)))/(x5 y5 z4 (x4 (y4 + y5 + y6) z3 + x5 (y5 + y6) (z3 + z4) + x6 y6 (z3 + z4 + z5)));
\end{lstlisting}

\begin{lstlisting}[language=Mathematica,caption={Example code}]
Simplify[K1 @@ arg1 == K9 @@ arg3]
(* True *)
\end{lstlisting}

\begin{lstlisting}[language=Mathematica,caption={Example code}]
Simplify[K3 @@ arg1 == K7 @@ arg3]
(* True *)
\end{lstlisting}

\begin{lstlisting}[language=Mathematica,caption={Example code}]
Simplify[K4 @@ arg1 == K6 @@ arg3]
(* True *)
\end{lstlisting}
\textbf{ Checking symmetries in F6: }

We can find the set of variables in a same way as what we did for $F_{9}^{(4)}$ and so here we omit most of the calculations.

Again as in $F_{9}^{(4)}$ , here we  specifically mean that the substitution 

\hspace*{-0.2cm} $\{x1, x2, x3, x4, x5, x6, y1, y2, y3, y4, y5 \} $ instead of $ \{ x6, x5, x4, x3, x2, x1,  $\\
$y5, y4 , y3, y2, y1 \}$ transforms $\tau_{6}^{(4)}$ to $\tau_{1}^{(4)}$, $\tau_{5}^{(4)}$ to $\tau_{2}^{(4)}$, $\tau_{4}^{(4)}$ to $\tau_{3}^{(4)}$ while leaving $F_{6}^{(4)}$ unchanged.

Therefore, those three pairs must enter the expression for $F_{6}^{(4)}$ symmetrically. 

\begin{lstlisting}[language=Mathematica,caption={Example code}]
arg1 = {a1, a2, a3, a4, a5, a6, b1, b2, b3, b4, b5};
arg2 = {a6, a5, a4, a3, a2, a1, b5, b4, b3, b2, b1};
\end{lstlisting}

\begin{lstlisting}[language=Mathematica,caption={Example code}]
F6[x1_, x2_, x3_, x4_, x5_, x6_, y1_, y2_, y3_, y4_, y5_] := (2 x1 x2 x5 x6 y2 (x2 y1 + x3 (y1 + y2)) y3^2 y4 (x5 y5 + x4 (y4 + y5)))/((x2 y2 + x1 (y1 + y2)) (x3 y3 + x2 (y2 + y3))^2 (x4 y3 + x5 (y3 + y4))^2 (x5 y4 + x6 (y4 + y5)));
\end{lstlisting}

\begin{lstlisting}[language=Mathematica,caption={Example code}]
F6 @@ arg1 == F6 @@ arg2
(* True *)
\end{lstlisting}

\begin{lstlisting}[language=Mathematica,caption={Example code}]
K1[x1_, x2_, x3_, x4_, x5_, x6_, y1_, y2_, y3_, y4_, y5_] := ((x2 y2 + x1 (y1 + y2)) (x3 y3 + x2 (y2 + y3)))/(x2 y2 (x3 y3 + x2 (y2 + y3) + x1 (y1 + y2 + y3)));
K2[x1_, x2_, x3_, x4_, x5_, x6_, y1_, y2_, y3_, y4_, y5_] := ((x2 y1 + x3 (y1 + y2)) (x3 y2 + x4 (y2 + y3)))/(x3 y2 (x2 y1 + (x3 + x4) (y1 + y2) + x4 y3));
K3[x1_, x2_, x3_, x4_, x5_, x6_, y1_, y2_, y3_, y4_, y5_] := ((x3 y3 + x2 (y2 + y3)) (x4 y4 + x3 (y3 + y4)))/(x3 y3 (x4 y4 + x3 (y3 + y4) + x2 (y2 + y3 + y4)));
K4[x1_, x2_, x3_, x4_, x5_, x6_, y1_, y2_, y3_, y4_, y5_] := ((x3 y2 + x4 (y2 + y3)) (x4 y3 + x5 (y3 + y4)))/(x4 y3 (x3 y2 + (x4 + x5) (y2 + y3) + x5 y4));
K5[x1_, x2_, x3_, x4_, x5_, x6_, y1_, y2_, y3_, y4_, y5_] := ((x4 y4 + x3 (y3 + y4)) (x5 y5 + x4 (y4 + y5)))/(x4 y4 (x5 y5 + x4 (y4 + y5) + x3 (y3 + y4 + y5)));
K6[x1_, x2_, x3_, x4_, x5_, x6_, y1_, y2_, y3_, y4_, y5_] := ((x4 y3 + x5 (y3 + y4)) (x5 y4 + x6 (y4 + y5)))/(x5 y4 (x4 y3 + (x5 + x6) (y3 + y4) + x6 y5));
\end{lstlisting}

\begin{lstlisting}[language=Mathematica,caption={Example code}]
Simplify[K1 @@ arg1 == K6 @@ arg2]
(* True *)
\end{lstlisting}

\begin{lstlisting}[language=Mathematica,caption={Example code}]
Simplify[K2 @@ arg1 == K5 @@ arg2]
(* True *)
\end{lstlisting}

\begin{lstlisting}[language=Mathematica,caption={Example code}]
Simplify[K3 @@ arg1 == K4 @@ arg2]
(* True *)
\end{lstlisting}

\footnote{Farrokh Razavinia, Department of discrete mathematics,  Moscow Institute of Physics and Technology, Institutskiy per., 9, Dolgoprudny, Moscow Oblast, Russia} \textsuperscript{1}

\footnote{Brendan B. Godfrey, Institute for Research in Electronics and Applied Physics (The University of Maryland), College Park, United States of America} \textsuperscript{2}

\end{document}